\documentclass[a4paper,12pt,reqno]{article}
\usepackage{amsmath,amsthm,amssymb,amscd}
\usepackage{amssymb}
\usepackage{amsthm}
\usepackage{amsmath}
\usepackage{latexsym}
\usepackage{epsfig}
\usepackage{amscd}
\usepackage{graphicx,wrapfig}
\usepackage{subcaption}
\usepackage[thinlines,thiklines]{easybmat}
\usepackage[table]{xcolor}
\usepackage{hyperref}
\usepackage{float}
\hypersetup{
    colorlinks=true,
    linkcolor=blue,
    filecolor=cyan,   
    urlcolor=magenta,
    citecolor=red,
    pdftitle={Overleaf Example},
    pdfpagemode=FullScreen,
}
\DeclareMathOperator{\lcm}{lcm}
\def\natu           {\mathbb N}
\def\inte 		{\mathbb Z}
\def\real		{\mathbb R}
\def\rati		{\mathbb Q}
\def\comp		{\mathbb C}

\title{Periodic Column Partial Sums\\ in the Riordan Array of a Polynomial}
 
\author{Nikolai A. Krylov \thanks{This work was supported by an AMS-Simons Research 
Enhancement Grant for Primarily Undergraduate Institution Faculty.}\\ ~ \\
Siena College, Department of Mathematics\\
515 Loudon Road, Loudonville NY 12211, USA\\ ~ \\
nkrylov@siena.edu}

\date{}

\begin{document}

\newtheorem{theorem}{Theorem}
\newtheorem{lemma}[theorem]{Lemma}
\newtheorem{claim}[theorem]{Claim}
\newtheorem{cor}[theorem]{Corollary}
\newtheorem{conj}[theorem]{Conjecture}
\newtheorem{prop}[theorem]{Proposition}
\newtheorem{question}{Problem}
\theoremstyle{definition}
\newtheorem{definition}[theorem]{Definition}
\newtheorem{example}[theorem]{Example}
\numberwithin{equation}{section}

\maketitle
\begin{abstract}
When $p(t)$ is a polynomial of degree $d$, $k$-th column of the Riordan array 
$\bigl(1/(1 - t^{d+1}), tp(t)\bigr)$ is an eventually periodic sequence with the 
repeating part beginning at the $1 + (k-1)(d+1)$-st term. The pre-periodic 
terms add up to the $(k-1)(d+1)$-st partial sum of the corresponding 
formal power series, and thus the Riordan array of $p(t)$ generates a 
sequence of column partial sums. We classify linear and quadratic polynomials, 
and present a particular family of polynomials of higher degrees, for which 
such sequences of column partial sums are eventually periodic.
\end{abstract}


\section{Introduction}\label{intro}

Let $p(t) = a_0+a_1 t + \cdots + a_d t^d \in \comp[t]$ be a polynomial of degree $d \geq 0$ 
with $a_0\neq 0$. Taking coefficients of the formal power series (f.p.s.)
$$
\frac{p(t)}{1- t^{d+1}} = \frac{a_0+a_1 t + \cdots + a_d t^d}{1- t^{d+1}} \in\comp[[t]]
$$
we create an infinite periodic sequence with repeating blocks 
$\{a_0,a_1,\ldots,a_d\}$. It is easy to see
that for every $k\in\natu_0 = \natu\cup\{0\}$, coefficients of the f.p.s.
\begin{equation}
\label{FPS}
\frac{\bigl(tp(t)\bigr)^k}{1 - t^{d+1}}.
\end{equation}
also make an eventually periodic sequence (Theorem 1, \cite{KrylovRA1}). We can 
put all these sequences together into a lower triangular matrix, where $k$-th column 
($k\geq 0$) is the coefficient sequence of (\ref{FPS}). Such a matrix is an example 
of a Riordan array, which is defined using a pair of f.p.s. $\bigr(f(t),g(t)\bigl)$. 
In our example above, 
$$
f(t) = \frac{1}{1-t^{d+1}},~~ \mbox{and} ~~ g(t) = tp(t).
$$
In general, a proper Riordan array over a field $\mathbb F$, is defined as a pair of 
two formal power series $\bigr(f(t),g(t)\bigl)$, where $f(t) = \sum_{i\geq 0} f_it^i\in \mathbb F[[t]]$ 
and  $g(t) \sum_{i\geq 0} g_it^i  \in \mathbb F[[t]] $. To guarantee nonzero elements along 
the main diagonal and the lower triangular form of the array, the order of $f(t)$ must be 0, 
that is $f_0 \neq 0$,  and the order of $g(t)$ must be 1, that is $g_0 = 0 \wedge g_1\neq 0$. 
In this paper we are interested in certain properties of the Riordan arrays 
\begin{equation}
\label{RAp}
\bigr(f(t),g(t)\bigl) = \left(\frac{1}{1-t^{d+1}}, tp(t)\right),
\end{equation}
where $p(t)$ is a polynomial of degree $d\geq 0$ with complex coefficients, 
and $p(0) \neq 0$. We call such an array as {\sl Riordan array of $p(t)$} and 
abbreviate it as $RAp$. For the general theory of Riordan arrays and the Riordan 
Group, we refer the reader to the original paper \cite{Shapiro1}, where this 
terminology was introduced, a thorough introductory text \cite{Barry}, and 
comprehensive monograph \cite{Shapiro2}. The latter one is a great resource 
on recent developments of the topic and the corresponding literature. For a brief 
introduction to the topic, we recommend survey articles \cite{Cameron2} and 
\cite{Davenport1}.

More specifically, I will continue here the study of Riordan arrays (\ref{RAp}) 
started in \cite{KrylovRA1}, and focus on polynomials $p(t)$, which generate 
circulant matrices of finite period (or finite order when the matrix is nonsingular). 
We follow standard terminology used in discrete dynamical systems. In particular, let 
$f : X \to X$ be a map, which gives the time evolution of the points of a set $X$. 
If $x\in X$, consider the iterates 
$$
x, ~ f(x), ~ f(f(x)), \ldots, 
$$
and denote the $n$-th iterate of $f$ at $x$ as $f^n(x)$, with convention 
that $f^0 = I$, is the identity map. ${\cal O}_ f^+(x)$ denotes the forward 
orbit of a point $x\in X$ under the iterates of $f$ , i.e.
$$
{\cal O}_ f^+(x) := \{x, f (x), f^2(x),\ldots, f^n(x),\ldots \} = \{f^n(x), n \in \natu_0\}.
$$
A point $x \in X$ is called {\sl preperiodic} with a period $\mu\in\natu$, if its 
orbit becomes periodic with period $\mu$ after a finite number of steps. That is, if 
$\exists k\in \natu_0$ such that 
$$
f^{\mu+j}(f^k(x)) = f^j(f^k(x)), ~ \forall j \in \natu_0.
$$
In such a case, we will also say that the orbit ${\cal O}_ f^+(x)$ is {\sl eventually periodic}. 
If $k = 0$, the point $x$ is called {\sl periodic}. Clearly, a point $x$ is preperiodic if 
and only if the cardinality of ${\cal O}_ f^+(x)$ is finite. If we have a polynomial 
$p(t) = a_0+a_1 t + \cdots + a_d t^d$ as above, 
we can write its coefficients as a row vector starting with $a_d$ 
$$
\nu = (a_d, a_{d-1}, \ldots ,a_0)\in \comp^{d+1},
$$
and consider a shift operator $T:\comp^{d+1}\to \comp^{d+1}$ defined as 
$$
T(a_d, a_{d-1},\ldots, a_1,a_0):= (a_0, a_d, a_{d-1},\ldots, a_1).
$$
Following \cite{Kra}, construct with this vector $\nu$ the $(d+1)\times (d+1)$ 
matrix $V_{p(t)} = {\rm circ} \{\nu\}$, 
and call it the {\sl circulant matrix associated to the polynomial $p(t)$}. 
Its rows are given by iterations of the shift operator acting on the vector $\nu$, i.e. the
$i$-th row of $V_{p(t)}$ is $T^i \nu,~ i\in\{0,\ldots,d \}$:
$$
V_{p(t)} = 
\begin{pmatrix}
a_d & a_{d-1} & \cdots & a_1 & a_0\\
a_0 & a_{d} & \cdots & a_2 & a_1\\
\vdots & \vdots & \ddots & \vdots & \vdots\\
a_{d-2} & a_{d-3} & \cdots & a_d & a_{d-1}\\
a_{d-1} & a_{d-2} & \cdots & a_0 & a_d\\
\end{pmatrix}.
$$
In particular, if $d=0$ and $p(t) = a_0$, then $V_{p(t)} = (a_0)$.
We will consider a discrete dynamical system $f : X \to X$, where 
$X = \comp^{d+1}$, the point $x = (0,0, \ldots,0,1)^T\in \comp^{d+1}$, 
and $f$ is the linear map given by $V_{p(t)}$ with 
\begin{equation}
\label{PerCirMat}
f^{\mu + j}(f^k(x)) = f^j(f^k(x)), ~~ \mbox{for some} ~ \mu > 1, ~ k\geq 0,
 ~~ \mbox{and} ~~ \forall j \in \natu_0.
\end{equation}
To simplify terminology, we will say that $V_{p(t)}$ is {\sl periodic with period 
$\mu$} if it satisfies (\ref{PerCirMat}), and $\mu$ is the smallest such positive integer.
Thus, we are interested in specific polynomials $p(t)\in\comp[t]$, which 
generate periodic circulant matrices $V_{p(t)}$. 
Take, for example, a polynomial $p(t) = 1/2 - t/2$. Then
$$
V_{p(t)} = 
\begin{pmatrix}
-1/2 & 1/2\\
1/2 & -1/2\\
\end{pmatrix}, ~ V^2_{p(t)} = 
\begin{pmatrix}
1/2 & -1/2\\
-1/2 & 1/2\\
\end{pmatrix}, ~ V^3_{p(t)} = V_{p(t)}.
$$
The orbit of $(0,1)^T$ equals 
$$
\left\{
\begin{pmatrix}
0 \\ 1 \end{pmatrix},~
\begin{pmatrix}
1/2 \\ -1/2 \end{pmatrix},~
\begin{pmatrix}
-1/2 \\ 1/2 \end{pmatrix}\right\}, 
$$
and both, the orbit and $V_{p(t)}$, have period 2.

According to Theorem 2 of \cite{KrylovRA1}, each column of  the Riordan array 
(\ref{RAp}) is an eventually periodic sequence. If we use the standard terminology, the top 
row of $RAp$ corresponds to $n = 0$, and the most left column is the 0-th column. Then 
the periodic cycle of $k$-th column ($k\geq 1$) begins at $1 + (k-1)(d+1)$-st place,  
and equals $k$-th iteration of $V_{p(t)}$ at $(0, \ldots,0,1)^T$ 
(see theorems 1 and 2 of \cite{KrylovRA1}). Here is the array corresponding to 
$p(t) = 1/2 - t/2$.
\begin{equation}
\label{matrixVPer}
\left(\frac{1}{1-t^2}, \frac{t-t^2}{2}\right) = 
\left(
\begin{BMAT}[2pt]{ccccccc}{cccccccc}
1 & 0 & 0 & 0 & 0 & 0 &  0  \\
0 & 1/2 & 0  & 0 & 0 & 0 &  0  \\
1 & -1/2 & 1/4 &  0 & 0 & 0 &  0  \\
0 & 1/2 & -1/2 & 1/8 & 0 & 0 &  0  \\
1 & -1/2 & 1/2 & -3/8 & 1/16 & 0 &  0  \\
0 & 1/2 & -1/2 & 1/2 & -1/4 & 1/32 & 0 \\
1 & -1/2 & 1/2 & -1/2 & 7/16 & -5/32 &  1/64  \\
\vdots & \vdots &  \vdots & \vdots & \vdots & \vdots & \ddots 
\addpath{(1,5,1)ruuldd}
\addpath{(2,3,1)ruuldd}
\addpath{(3,1,1)ruuldd}
\end{BMAT}
\right)
\end{equation}

Polynomials that generate periodic $V_{p(t)}$, produce {\sl doubly periodic} Riordan 
arrays, because in addition to the eventual periodicity in each column ({\sl vertical} 
periodicity), there is {\sl horizontal repetition} too, as one can see it in (\ref{matrixVPer}). 
If we look at the {\sl preperiodic} part in each column of such Riordan array, we 
notice that for some polynomials, the {\sl preperiodic} terms add up to numbers, 
which behave periodically as well. In particular, returning to our example of $p(t) = 1/2 - t/2$, 
if we add first $2k - 1$ terms (starting in the 0-th row) in $k$-st column, we'll 
get the following sequence of column partial sums
$$
\left\{0,~\frac{1}{4},~\frac{-1}{4},~\frac{1}{4},~\frac{-1}{4},~\frac{1}{4},~\frac{-1}{4},\ldots\right\}.
$$
According to the Fundamental Theorem of Riordan Arrays, (the FTRA, see \S 5.1 
of \cite{Barry} or Theorem 3.1 of \cite{Shapiro2}), we can think of the $k$-th 
column of $RAp$ as a sequence of coefficients of the formal power series 
\begin{equation}
\label{CFPS}
\sum\limits_{i=0}^{\infty} C_{i,k} t^i  = \frac{\bigl(tp(t)\bigr)^k}{1 - t^{d+1}}.
\end{equation}
Then, if we add all elements in the {\sl preperiodic} part of $k$-th column 
(starting with 0 in the 0-th row), we will get exactly the $(k-1)(d+1)$-st partial 
sum of this infinite series
\begin{equation}
\label{PSum}
S_{(k-1)(d+1)} = \sum\limits_{i=0}^{(k-1)(d+1)} C_{i,k}.
\end{equation}

\begin{definition}
Given the Riordan array (\ref{RAp}) of a degree $d$ polynomial $p(t)$, the 
$(k-1)(d+1)$-st partial sum (\ref{PSum}) of the formal power series (\ref{CFPS})
will be denoted by $S_{[k]}$ and called the {\sl partial sum of $k$-th column}, 
or simply as the {\sl column partial sum}. 
\end{definition}

Here are two more examples from \cite{KrylovRA1}, where I considered polynomials 
with real coefficients. If $p(t) = (-1 + 2t + 2t^2)/3$, then the sequence of column 
partial sums $\left\{S_{[k]}\right\}_{k\geq 1}$ has period 6 (see Proposition 5 
in \cite{KrylovRA1}, or Proposition \ref{PropD} below) with the repeating part
$$
\left\{0,~\frac{-1}{3},~\frac{-2}{3},~\frac{-2}{3},~\frac{-1}{3},~ 0 \right\}.
$$
If $p(t) = (2 - t + 2t^2)/3$, the sequence of partial sums is not eventually periodic, 
and the first few terms of the sequence are as follows.
$$
\left\{S_{[k]}\right\}_{k\geq 1} = 
\left\{0,~0,~\frac{1}{3},~1,~\frac{5}{3},~2, ~2,~2, ~\frac{7}{3},~3,\cdots \right\}
$$

The main objective of this paper is to find polynomials $p(t)\in \comp[t]$, 
which generate eventually periodic sequences of the column partial sums 
\begin{equation}
\label{SeqPSums}
\left\{S_{[k]}\right\}_{k\geq 1}
\end{equation}
in the Riordan array (\ref{RAp}). Linear polynomials with this property 
are given in Theorem \ref{theorem10}, 
and a particular family of higher degree polynomials of such kind is presented 
in Proposition \ref{PropD}. Theorem \ref{thm2} below classifies the corresponding 
quadratic polynomials (for the proofs see Propositions \ref{propA}, \ref{propB}, 
and \ref{propC} in \S 3).

\begin{theorem}
\label{thm2}
Let $p(t) = a + bt + ct^2$ be a polynomial over $\comp$ with $a\cdot c\neq 0$. Denote 
the primitive root of unity $e^{2\pi i/3}$ by $\xi$, and the eigenvalues of $V_{p(t)}$ 
by $\lambda_i$, where 
$$
\lambda_1 = a + b + c, ~~~ \lambda_2 = a\xi^2 + b\xi + c, ~~ \mbox{and} 
~~ \lambda_3 = a\xi + b\xi^2 + c.
$$
\begin{itemize}
\item If $\lambda_1\neq 0$, then $V_{p(t)}$ has a finite period $\mu$ and the 
column partial sums are eventually periodic if and only if 
$$
p(t) = a(1 - 2t - 2t^2), ~~ (3a)^{\mu} = 1,~~  \mbox{and} ~~ 6~|~\mu. 
$$
\item If $\lambda_1 = \lambda_2 = 0$, then $V_{p(t)}$ has a finite period $\mu$ 
and the sequence of column partial sums is eventually periodic with the 
general term 
$$
S_{[k]} = (3a\xi)^k\frac{(\xi - 1)}{9\xi},~ k\geq 2
$$ 
if and only if 
$$
p(t) = a(1 - t)(1 - \xi t) ~~ \mbox{and} ~~ (3a\xi)^{\mu} = 1.
$$
\item If $\lambda_1 = \lambda_3 = 0$, then $V_{p(t)}$ has a finite period $\mu$ 
and the sequence of column partial sums is eventually periodic with the 
general term 
$$
S_{[k]} = (3a\xi^2)^k\frac{(1 - \xi)}{9}, ~ k\geq 2
$$
if and only if 
$$
p(t) = a(1 - t)(\xi^2 - \xi t) ~~ \mbox{and} ~~ (3a\xi^2)^{\mu} = 1.
$$
\item If $\lambda_1 = 0$ and $\lambda_2\cdot\lambda_3 \neq 0$, then 
$V_{p(t)}$ has a finite period $\mu$ and the column partial sums 
are eventually periodic if and only if 
$$
p(t) = a(1 - t)\bigl(1 + (r+1)t\bigr),
$$
where $a$ and $r$ satisfy the system
$$
\left\{
\begin{array}{lcl}
\bigl((\xi - 1)a(r - \xi^2)\bigr)^{\mu} & = & 1\\ 
\bigl((\xi - 1)a(1 - \xi^2r)\bigr)^{\mu} & = & 1.\\ 
\end{array}
\right.
$$ 
In this case, $r\in\real$ is a root of a polynomial of degree $\mu$ with rational 
coefficients.
\end{itemize}
\end{theorem}

\noindent It is interesting to note, that 1 is a root of all these polynomials $p(t)$, with
$\lambda_1 = 0$. If $\lambda_1 \neq 0$, the roots of $p(t)$ are irrational numbers 
$(-1\pm\sqrt{3})/2$, $\forall a \in\comp\setminus \{0\}$.

~

The paper is organized as follows. Section 2, contains several technical lemmas, 
where in particular, we show how to identify the column partial sums as coefficients 
of certain formal power series. In Section 3, we discuss the main results for 
linear and quadratic polynomials: Theorem \ref{theorem10}, and Propositions 
\ref{propA}, \ref{propB}, and \ref{propC} (i.e. Theorem \ref{thm2}). Section 4 shows 
that the column partial sums of the polynomial 
$p(t) = a\bigl((d - 1) - 2t - 2t^2 - \cdots - 2t^d\bigr)$ are eventually periodic, 
when $d\geq 2$ and the product $a\cdot (d+1)$ is a root of unity. 
In the last section, we visualize several examples of eventually periodic 
column partial sums as graphs in the plane. 

~

Throughout the paper, if $p(t)$ has degree $d$, then $(d+1)$-st primitive root of unity 
$e^{2\pi i/(d+1)}$ will be denoted by $\xi$, unless specifically stated otherwise.

%%%%%%%%%%%%%%%%%%%%%%%%%%%
%%%%%%%%%%%%%%%%%%%%%%%%%%%
%%%%%%%%                       %%%%%%%%%%%%
%%%%%%%%     Section 2   %%%%%%%%%%%%
%%%%%%%%                       %%%%%%%%%%%%
%%%%%%%%%%%%%%%%%%%%%%%%%%%
%%%%%%%%%%%%%%%%%%%%%%%%%%%

\section{Riordan array of column partial sums and auxiliary lemmas}

First of all, for a constant non-zero polynomial $p(t) = a$ the sequence of partial sums 
(\ref{SeqPSums}) consists of zeroes only. This is true for all $a\in \comp\setminus\{0\}$.
$$
RAa = \left( \frac{1}{1-t}, ta\right) = 
\begin{pmatrix}
1 & 0 & 0 & 0 & 0 & \cdots \\
1 & a & 0 & 0 & 0 & \cdots \\
1 & a & a^2 & 0 & 0 & \cdots \\
1 & a & a^2 & a^3 & 0 & \cdots \\
%1 & a & a^2 & a^3 & a^4 & \cdots \\
 \vdots & \vdots & \vdots & \vdots & \vdots & \ddots \\
\end{pmatrix}
$$
Thus, from now on assume that $p(t)$ has degree $d\geq 1$, and denote the Riordan array 
(\ref{RAp}) of $p(t)$ by $(f,g)$. If we consider all partial sums $\{S_n\}_{n\geq 0} ^{\infty}$ 
in each column of $RAp$, we will obtain another Riordan array, which equals the product
\begin{equation}
\label{product1}
\left(\frac{1}{1-t},t\right) \cdot \left(\frac{1}{1-t^{d+1}}, tp(t)\right) = 
\left(\frac{1}{(1-t)(1-t^{d+1})}, tp(t)\right).
\end{equation}
This equality follows from an observation that the $n$-th partial sum $S_n$ of the 
$m$-th column can be presented as a product of the infinite row vector 
$\{1,\ldots,1,0,0,\ldots\}$, which has 1 in the fist $n$ positions and 0 after that,
with the $m$-th column of $RAp$ (cf.  Proposition 3 of \cite{Barry2}, or  
Theorem 2.7 of \cite{Shapiro2}). Hence our sequence 
of partial sums (\ref{SeqPSums}) is a sequence of particular elements of the 
Riordan array (\ref{product1}), and we can apply the ``{\sl coefficient of} " operator to 
obtain some formulas for the elements of (\ref{SeqPSums}). To be precise, since we 
are adding the first $(k-1)(d+1)+1$ elements in the $k$-th column, we have 
$$
S_{[k]} = S_{(k-1)(d+1)} = \sum\limits_{i=0}^{(k-1)(d+1)} C_{i,k}.
$$
If we denote entries of the matrix product (\ref{product1}) by $h_{n,k}$, then we have
$$
\left\{S_{[k]}\right\}_{k\geq 1} = \{h_{0,1},~h_{d+1,2},~ h_{2(d+1),3}, ~ 
\ldots, ~ h_{k(d+1),k+1}, ~\ldots \} = \left\{h_{l(d+1),l+1}\right\}_{l\geq 0}.
$$
Thus, using the ``{\sl coefficient of} " operator we can write 
\begin{equation}
\label{ParSumCoef}
S_{[k]} = h_{(k-1)(d+1),k} = [t^{(k-1)(d+1)}]\frac{(tp(t))^k}{(1-t)(1-t^{d+1})},~~ k\geq 1.
\end{equation}

\noindent Note that $S_{[1]} = 0$ for any $p(t)$ with $a_0\neq 0$. It is well-known 
that rational generating functions give linear recurrence sequences, and vice versa 
(see \cite{Lando}, \S 2.3). The following simple result characterizes eventually 
periodic sequences. 

\begin{lemma}
\label{lemma4}
Generating function of a sequence $\{s_i\}_{i\geq 0}$ can be written as a rational function
\begin{equation}
\label{RatGF}
\sum\limits_{i=0}^{\infty} s_it^i = \frac{Q(t)}{1-t^n}
\end{equation}
for some polynomial $Q(t)\in \comp[t]$
if and only if the sequence $\{s_i\}_{i\geq 0}$ is eventually periodic with the period $n$.
\end{lemma}
\begin{proof}
Suppose the periodicity starts at the term $s_k$, and has the repeating part  
$\{s_k,\ldots, s_{k + n - 1}\}$ of length $n$. Consider two polynomials
$$
q_0(t) := s_{0} +s_{1}t + \cdots + s_{k-1}t^{k-1}, ~ \mbox{and} ~ 
q(t):= s_k +s_{k+1}t + \cdots + s_{k + n - 1}t^{n-1},
$$
then
$$
\sum\limits_{i=0}^{\infty} s_it^i = q_0(t) + q(t)t^{k} + q(t)t^{k + n} + q(t)t^{k + 2n} + \cdots = 
\frac{q_0(t)(1 - t^n) + q(t)t^k}{1 - t^n}.
$$

On the other hand, assume that (\ref{RatGF}) holds true. If we write the polynomial $Q(t)$ as
$q_0 + q_1t + \cdots + q_kt^k$, then 
\begin{equation}
\label{EQ2}
\sum\limits_{i=0}^{\infty} s_it^i = \frac{Q(t)}{1-t^n} = \sum\limits_{a=0}^{\infty} 
(q_0 + q_1 + \cdots + q_kt^k)t^{na},
\end{equation}
and if $k < n$, the sequence $\{s_i\}_{i\geq 0}$ is clearly periodic, with the repeating part  
$$
\{q_0,q_1,\ldots, q_k, \underbrace{0,\ldots,0}_{n-k-1}\}.
$$
If $k\geq n$, then let $k = pn + r,~0\leq r < n,$ and present $Q(t)$ as
$$
Q(t) = \sum\limits_{i=0}^{p-1}Q_{i}(t)t^{in} + R(t)t^{pn},
$$
where 
$$
Q_i(t) = q_{in}+q_{in+1}t + \cdots + q_{(i+1)n-1}t^{n-1},
$$
and 
$$
R(t) = q_{pn}+q_{pn+1}t + \cdots + q_{pn + r}t^{r}.
$$
Then we can rewrite (\ref{EQ2}) as
$$
\sum\limits_{i=0}^{\infty} s_it^i = Q_0(t) + \Bigl(Q_0(t) + Q_1(t)\Bigr)t^n + 
\Bigl(Q_0(t) + Q_1(t) + Q_2(t)\Bigr)t^{2n} + 
$$
$$
\cdots + \Bigl(Q_0(t) + \cdots + Q_{p-1}(t)\Bigr)t^{(p-1)n} + 
\Bigl(Q_0(t) + \cdots + Q_{p-1}(t) + R(t)\Bigr)\frac{t^{pn}}{1-t^n}
$$
to see that the sequence $\{s_i\}_{i\geq 0}$ is indeed eventually periodic. The repeating
part here is a cyclic permutation of the coefficients of the polynomial
$$
\Pi(t) : = Q_0(t) + \cdots + Q_{p-1}(t) + R(t),
$$
and if we write this polynomial as 
$$
\Pi(t) = p_0 + p_1t + \cdots + p_{n-1}t^{n-1},
$$
the repeating part will equal $\{p_{r+1}, \ldots, p_{n-1},p_0,\ldots p_r\}$.
\end{proof}

Notice, that $Q(t)$ and $1-t^n$ may have common divisors. 
For example, if $\xi = e^{2\pi i/3}$, the generating 
function for the periodic sequence with the repeating part $\{1,\xi,\xi^2\}$ will be 
$$
\sum\limits_{i=0}^{\infty} s_it^i = \frac{1+t\xi + t^2\xi^2}{1-t^3} = 
\frac{\xi^2(1 - t)(\xi - t)}{(1 - t)(\xi - t)(\xi^2 - t)} = \frac{\xi^2}{\xi^2 - t} = \frac{1}{1 - \xi t}.
$$

Given a polynomial $p(t) = a_0 + \cdots + a_dt^d$, let us associate 
with each power of the matrix $V_{p(t)}$ a particular polynomial as follows.

\begin{definition}
\label{definition3}
Let $p(t) = a_0+a_1 t + \cdots + a_d t^d \in \comp[t]$ be a polynomial of degree 
$d\geq 0$ with $a_0\neq 0$, and $V_{p(t)}$ its associated circulant matrix. 
For each $k\in\natu_0$ define the polynomial $p_k(t)$ by the following formula.
$$
p_k(t) = \sum\limits_{i = 0}^d a_{k,i}t^i : = (1,t,\ldots, t^d) \cdot V^k_{p(t)} \cdot \begin{pmatrix}
a_0\\
a_1\\
\vdots \\
a_d\\
\end{pmatrix}
$$
\end{definition}

Next we state and prove several technical lemmas, which will be used in the following 
sections. Until section 4, $\xi$ will stand for $e^{2\pi i/3}$.

\begin{lemma}
\label{lemma2}
$$
p_k(1) = \left( \sum\limits_{i = 0}^d a_i\right)^{k+1},~\forall k\in\natu_0.
$$
\end{lemma}
\begin{proof}
The base of induction with $k=0$ is clear since $p_0(t) = p(t)$. Assume now that the 
statement holds true for $k$. Then 
$$
p_{k+1}(1) = (1,1,\ldots,1) \cdot V^{k+1}_{p(t)} \cdot 
\begin{pmatrix}
a_0\\
a_1\\
\vdots \\
a_d\\
\end{pmatrix} = \bigl((1,1,\ldots,1) \cdot V_{p(t)} \bigr)\cdot 
\begin{pmatrix}
a_{k,0}\\
a_{k,1}\\
\vdots \\
a_{k,d}\\
\end{pmatrix}
$$
$$
= \left(\left( \sum\limits_{i = 0}^d a_i\right)(1,1,\ldots,1)\right)\cdot 
\begin{pmatrix}
a_{k,0}\\
a_{k,1}\\
\vdots \\
a_{k,d}\\
\end{pmatrix} = \left( \sum\limits_{i = 0}^d a_i\right)(a_{k,0} + a_{k,1} + \cdots + a_{k,d})
$$
$$
= \left( \sum\limits_{i = 0}^d a_i\right)\cdot \left( \sum\limits_{i = 0}^d a_i\right)^{k+1} = 
\left( \sum\limits_{i = 0}^d a_i\right)^{k+2}.
$$
\end{proof}

\begin{lemma}
\label{lemma5}
Suppose we have two nonzero complex numbers 
$z$ and $w$ s.t.
$$
|w| = 1, ~~ \mbox{and} ~~ | z - w | = | z - \xi w |.
$$
Then $z$ lies on the line through the origin, bisecting 
the angle $\theta = 2\pi/3$ between the vectors representing $w$ and $\xi w$.
\end{lemma}
\begin{proof}
Represent the numbers $z, w$, and $\xi w$ graphically in the complex plane by the 
points $A$, $B_1$, and $B_2$ correspondingly. Denote the origin by $O$, and consider 
two triangles $\Delta OAB_1$ and $\Delta OAB_2$. Let $\alpha_i$, $i\in\{1,2\}$ 
denote the angle $\alpha_i = \angle AOB_i$. Then our 
assumptions together with the law of cosines imply that
$\cos \alpha_1 = \cos \alpha_2$, and since $\alpha_i < \pi$, we see that 
$\alpha_1 = \alpha_2$.
\end{proof}

\begin{lemma}
\label{lemma6}
The following identity holds true for all $k \geq 2$.
\begin{equation}
\label{EQ9}
[t^{2k - 3}]\frac{(1-t^2)^{k - 2}(1+t)^2}{(1+t+t^2)} = 
2(-\sqrt{3})^{k - 3} \sin\left(\frac{(k - 4)\pi}{6}\right). 
\end{equation}
\end{lemma}
\begin{proof}
We will use the binomial theorem for $(1-t^2)^{k-2}$, together with the 
power series expansion for $(1+t)^2/(1+t+t^2)$, so let us write
$$
[t^{2k - 3}]\frac{(1-t^2)^{k - 2}(1+t)^2}{(1+t+t^2)}= 
[t^{2(k - 1)}]\left((1-t^2)^{k - 2}\cdot \frac{t(1+t)^2}{(1+t+t^2)}\right).
$$
Since we need the coefficient of the even power of $t$ and the 
binomial expansion will have only the even powers of $t$, we can ignore the 
odd-power terms in the expansion of 
$$
\frac{t(1+t)^2}{(1+t+t^2)} = t + \frac{t^2(1-t)}{1- t^3} = 
t + t^2(1-t)\cdot(1 + t^3 + t^6 + \cdots) 
$$
$$
= t + t^2 - t^3 + t^5 - t^6 + t^8 - t^9 + t^{11} - t^{12} + \cdots .
$$
Therefore, the L.H.S. of (\ref{EQ9}) equals
\begin{equation}
\label{EQ9b}
[t^{2(k - 1)}]\sum\limits_{i=0}^{k-2} (-1)^i{k-2\choose i} t^{2i}\cdot
\bigl(t^2 - t^6 + t^8 - t^{12} + t^{14} - t^{18} + \cdots\bigr).
\end{equation}
Replacing $t^2$ with $x$ in (\ref{EQ9b}) and using the convolution, we rewrite it as 
$$
[x^{k - 1}]\sum\limits_{i=0}^{k-2} (-1)^i{k-2\choose i} x^i \cdot
\bigl(x - x^3 + x^4 - x^6 + x^7 - x^9 + \cdots\bigr) 
$$
\begin{equation}
\label{EQ10}
= \sum\limits_{n = 0}^{k - 2} (-1)^n {k-2\choose n}\cdot [x^{k-1-n}]
\bigl(x - x^3 + x^4 - x^6 + x^7 - x^9 + \cdots\bigr).
\end{equation}
The coefficients of the f.p.s. in (\ref{EQ10}) make a periodic sequence with the 
repeating blocks $\{1,0,-1\}$, and the imaginary parts of the odd powers of 
$\eta = e^{\frac{2\pi i}{6}}$ make a periodic sequence with the repeating block
$$
\left\{{\rm Im}\bigl(\eta), {\rm Im}\bigl(\eta^3),{\rm Im}\bigl(\eta^5)\right\} = 
\left\{\frac{\sqrt{3}}{2}, 0, \frac{-\sqrt{3}}{2} \right\}.
$$
Hence we can write the f.p.s. expansion in (\ref{EQ10}) as
$$
\frac{x(1 + x)}{1 + x + x^2} = x - x^3 + x^4 - x^6 + x^7 - x^9 + \cdots = 
\sum\limits_{n=1}^{\infty}\frac{2}{\sqrt{3}} {\rm Im}(\eta^{2n-1}) x^n,
$$
and deduce that the L.H.S. of (\ref{EQ9}) equals
$$
\sum\limits_{n = 0}^{k - 2} (-1)^n {k-2\choose n}\cdot 
\frac{2}{\sqrt{3}} {\rm Im}\bigl(\eta^{2(k-1-n)-1}\bigr)
$$
\begin{equation}
\label{EQ11}
= {\rm Im}\left(\frac{2\eta}{\sqrt{3}}  \sum\limits_{n = 0}^{k - 2} {k - 2\choose n}(-1)^n \cdot 
\xi ^{(k-2-n)}\right) = 
{\rm Im}\left(\frac{2\eta}{\sqrt{3}} (\xi - 1)^{k-2}\right).
\end{equation}
Applying de Moivre's formula to $\xi - 1 = -\sqrt{3}e^{\frac{-\pi i}{6}}$, we obtain 
$$
\eta (\xi - 1)^{k-2} = e^{\frac{2\pi i}{6}} (-\sqrt{3})^{k-2} e^{\frac{(2-k)\pi i}{6} } = 
(-\sqrt{3})^{k-2} e^{\frac{( 4 - k)\pi i}{6} },
$$
which implies that (\ref{EQ11}) equals
$$
2 (-\sqrt{3})^{k - 3} \sin\left(\frac{(k - 4)\pi}{6}\right),
$$
and finishes the proof of (\ref{EQ9}).
\end{proof}

\begin{lemma}
\label{lemma7}
The following identities hold true for all $k \geq 2$ and $i \in \{0,1, 2, 3\}$.
\begin{equation}
\label{EQ3}
\left[t^{2k-i}\right] \frac{(1-t)^{k-2}(1 - \xi t)^{k-1}}{1- \xi^2t} = 3^{k-2}(\xi -1 ) \xi^{k + i - 1},
\end{equation}
and 
\begin{equation}
\label{EQ3B}
\left[t^{2k-i}\right] \frac{(1-t)^{k-2}(1 - \xi^2 t)^{k-1}}{1- \xi t} = 3^{k-2}(1 - \xi) \xi^{2k - i},
\end{equation}
\end{lemma}
\begin{proof}
We will prove the first one by the induction on $k$, and leave the verification of 
(\ref{EQ3B}) to the reader. Since $1 - t^3 = (1-t)(1-\xi t)(1 - \xi^2 t)$, according 
to Lemma \ref{lemma4}, the rational function
$$
\frac{1 - \xi t}{1 - \xi^2 t} = \frac{(1 - t)(1 - \xi t)^2}{1 - t^3} = 
\frac{1 + (\xi^2 - \xi)t + (\xi - 1)t^2  - \xi^2t^3}{1 - t^3} 
$$ 
\begin{equation}
\label{EQ4}
= 1 + \sum\limits_{i = 0} ^{\infty} \bigl((\xi^2-\xi)t + (\xi - 1)t^2 + (1 - \xi^2)t^3\bigr)t^{3i}
\end{equation}
is a generating function of an eventually periodic sequence with the repeating part  
$\{\xi^2 - \xi, \xi - 1, 1 - \xi^2\}$. Therefore the base of the induction follows directly 
from (\ref{EQ4}). Assume now that the equality (\ref{EQ3}) holds true for an 
arbitrary $k\geq 2$ and $0\leq i\leq 3$, and consider
$$
\left[t^{2(k + 1) - i}\right] \frac{(1-t)^{k-1}(1 - \xi t)^{k}}{1- \xi^2t} = 
\left[t^{2k + 2 - i}\right] \frac{(1-t)^{k-2}(1 - \xi t)^{k - 1}}{1- \xi^2t} (1 + \xi^2 t + \xi t^2).
$$
If $i\geq 2$, let $j:=i - 2\in \{0,1\}$ and use the induction hypothesis together with the 
convolution rule to obtain
$$
\left[t^{2k - j}\right] \frac{(1-t)^{k-2}(1 - \xi t)^{k - 1}}{1- \xi^2t} (1 + \xi^2 t + \xi t^2) 
$$
$$
= 3^{k-2}(\xi -1)\xi^{k-1+j} + 3^{k-2}(\xi -1)\xi^{k+j} \xi^2 + 3^{k-2}(\xi -1)\xi^{k+1+j} \xi 
$$
$$
= 3^{k-2}(\xi -1)\xi^{k-1+j}(1 + \xi^3 + \xi^3) = 3^{k-1}(\xi -1)\xi^{k + i - 3} = 3^{k-1}(\xi -1)\xi^{k + i}.
$$
This proves (\ref{EQ3}) for $i\in \{2,3\}$. Next, for $i=1$ we have 
$$
\left[t^{2k + 1}\right] \frac{(1-t)^{k-1}(1 - \xi t)^{k}}{1- \xi^2t} = 
\left[t^{2k + 1}\right] \frac{(1-t)^{k-2}(1 - \xi t)^{k - 1}}{1- \xi^2t} (1 + \xi^2 t + \xi t^2) 
$$
\begin{equation}
\label{EQ5}
= \left[t^{2k + 1}\right] \frac{(1-t)^{k-2}(1 - \xi t)^{k - 1}}{1- \xi^2t} + 
3^{k-2}(\xi -1)\xi^{k - 1}\xi^2 + 3^{k-2}(\xi -1)\xi^{k} \xi.
\end{equation}
Since the rational function 
\begin{equation}
\label{EQ6}
\frac{(1-t)^{k-2}(1 - \xi t)^{k - 1}}{1- \xi^2t} = \frac{(1-t)^{k-1}(1 - \xi t)^{k}}{1 -t^3}
\end{equation}
generates an eventually periodic sequence with period 3, and the repeating part 
begins with the coefficient of $t^{2k-3}$ (recall the end of proof of Lemma \ref{lemma4}), 
$$
\left[t^{2k + 1}\right] \frac{(1-t)^{k-2}(1 - \xi t)^{k - 1}}{1- \xi^2t} = 
\left[t^{2k - 2}\right] \frac{(1-t)^{k-2}(1 - \xi t)^{k - 1}}{1- \xi^2t} = 3^{k-2}(\xi -1)\xi^{k + 1}.
$$
Thus, (\ref{EQ5}) equals
$$
3^{k-2}(\xi -1)\xi^{k + 1} + 3^{k-2}(\xi -1)\xi^{k - 1}\xi^2 + 3^{k - 2}(\xi -1)\xi^{k} \xi = 
3^{k-1}(\xi -1)\xi^{k+ 1},
$$
as required to prove the case when $i = 1$. The final case of $i = 0$ is proved in 
exactly the same way using the induction hypothesis together with Lemma \ref{lemma4}  
and the periodicity of the rational function in 
(\ref{EQ6}).
\end{proof}

Notice that (\ref{EQ3}) is not in general true when $i > 3$. For example, 
for $k = 3$ and $i = 4$ we have
$$
[t^2]\frac{(1-t)(1 - \xi t)^2}{(1 - \xi^2t)} = 3\xi - 2 \neq 3(\xi - 1).
$$

\begin{lemma}
\label{lemma8}
For any fixed $m\in\natu$, 
\begin{equation}
\label{difBIN}
(x\xi - 1)^m - (\xi - x)^m = z_m\cdot f_m(x),
\end{equation}
where $\xi= e^{2\pi i/3}$, $z_m\in\comp\setminus\{0\}$ and $f_m(x)\in\rati[x]$.
\end{lemma}
\begin{proof}
Using the binomial theorem write the difference in (\ref{difBIN}) as
$$
\sum\limits_{i=0}^m {m\choose i} \bigl( (x\xi)^i(-1)^{m-i} - (-x)^i \xi^{m-i}\bigr) = 
\sum\limits_{i=0}^m (-x)^i {m\choose i} \bigl( \xi^i(-1)^{m} - \xi^{m-i}\bigr).
$$
It follows easily from $1 + \xi + \xi^2 = 0$, that 
$$
(-1)^m - \xi^m = 0 \Longleftrightarrow 2 - (-1)^m\bigl(\xi^m + \xi^{2m}\bigr) = 0 
\Longleftrightarrow 6 ~ | ~ m,
$$
and so the sequence
$$
\left\{2 - (-1)^m\bigl(\xi^m + \xi^{2m}\bigr)\right\}_{m\geq 0} = 
\{0, 1, 3, 4, 3, 1, 0, 1, 3, 4, 3, 1, 0, \ldots\}
$$
is periodic with period 6. If 6 doesn't divide $m$, denote $(-1)^m - \xi^m$ by $z$, 
then $\bar{z} = (-1)^m - \xi^{2m}$, $z\bar{z} = 2 - (-1)^m(\xi^m + \xi^{2m})\in\inte^+$, 
and we can write
$$
\xi^i(-1)^{m} - \xi^{m-i} = z\frac{(\xi^i(-1)^{m} - \xi^{m-i})\bar{z}}{z\bar{z}} = 
z \frac{(\xi^i(-1)^{m} - \xi^{m-i})\cdot((-1)^m - \xi^{2m})}{2 - (-1)^m\bigl(\xi^m + \xi^{2m}\bigr)} 
$$
$$
= z\frac{\xi^i - (-1)^m\xi^{m-i} - (-1)^m\xi^{2m+i} + \xi^{-i}}{2 - (-1)^m\bigl(\xi^m + \xi^{2m}\bigr)} 
$$
\begin{equation}
\label{ZCoeff}
= z\frac{\left(\xi^i + \bar{\xi}^i\right) - (-1)^m\left(\xi^{m-i} +\bar{\xi}^{m - i}\right)}
{2 - (-1)^m\bigl(\xi^m + \bar{\xi}^m\bigr)} = z\cdot q_{mi},
\end{equation}
where $q_{mi}\in\rati$ is a rational number, a priori depending on $m$ and $i$. When 
$6 ~ | ~ m$, we have $\xi^i(-1)^{m} - \xi^{m-i} = \xi^i  - \bar{\xi}^i = \sqrt{-3}\cdot q_{6i}$, 
where $q_{6i}\in\{-1,0,1\}$. Therefore, we can write the difference in (\ref{difBIN}) as 
$$
z\sum\limits_{i=0}^m (-x)^i {m\choose i} q_{mi},
$$
and use 
$$
z_m : = 
\left\{
\begin{array}{rcl}
(-1)^m - \xi^m & \mbox{if} & 6 \nmid m \\
\sqrt{-3} & \mbox{if} &6 ~ | ~m\\
\end{array}
\right.
$$
and 
$$
f_m(x) : = \sum\limits_{i=0}^m (-x)^i {m\choose i} q_{mi} \in \rati[x],
$$
where $q_{mi} = q_{6i}$ if $6~|~m$. 
\end{proof}

%%%%%%%%%%%%%%%%%%%%%%%%%%%
%%%%%%%%%%%%%%%%%%%%%%%%%%%
%%%%%%%%                       %%%%%%%%%%%%
%%%%%%%%     Section 3   %%%%%%%%%%%%
%%%%%%%%                       %%%%%%%%%%%%
%%%%%%%%%%%%%%%%%%%%%%%%%%%
%%%%%%%%%%%%%%%%%%%%%%%%%%%

\section{Linear and Quadratic polynomials}

\subsection{Linear polynomials}

According to Theorem 4. of \cite{KrylovRA1}, if $p(t) = a + bt$, then
$$
V^{\mu}_{p(t)}
\begin{pmatrix} 
a \\ b \\
\end{pmatrix} = 
\begin{pmatrix}
b & a\\
a & b\\
\end{pmatrix}^{\mu} 
\begin{pmatrix} 
a \\ b \\
\end{pmatrix}  = \frac{1}{2}
\begin{pmatrix}
1 & -1\\
1 & 1
\end{pmatrix}
\begin{pmatrix}
(a+b)^{\mu+1}\\
(b-a)^{\mu+1}
\end{pmatrix},
$$
where $a+b$ and $b-a$ are the eigenvalues of $V_{p(t)}$. If $V_{p(t)}$ has period $\mu$, 
equation (\ref{PerCirMat}) implies 
\begin{equation}
\label{PerCirMatLin}
\frac{1}{2}
\begin{pmatrix}
1 & -1\\
1 & 1
\end{pmatrix}
\begin{pmatrix}
(a+b)^{\mu +1}\\
(b-a)^{\mu +1}
\end{pmatrix} = 
\begin{pmatrix}
a\\
b\\
\end{pmatrix} \Longleftrightarrow 
\begin{pmatrix}
(a+b)^{\mu +1}\\
(b-a)^{\mu +1}
\end{pmatrix} = 
\begin{pmatrix}
a+b\\
b-a\\
\end{pmatrix},
\end{equation}
so we separate three cases (recall that $ab\neq 0$):
\begin{itemize}
\item[]{Case (1)} If $a = b$ then 
$$
(2a)^{\mu} = 1 \Longrightarrow a = b = \frac{1}{2}e^{\frac{2\pi i s}{\mu}},~s\in\{1,\ldots,\mu\}.
$$
\item[]{Case (2)} If $a = - b$ then 
$$
(2b)^{\mu} = 1 \Longrightarrow - a = b = \frac{1}{2}e^{\frac{2\pi i s}{\mu}},~s\in\{1,\ldots,\mu\}.
$$
\item[]{Case (3)} If none of the eigenvalues is zero, then $(a+b)^{\mu} = (b-a)^{\mu}  = 1$ and 
$$
(a,b) = \left(\frac{e^{\frac{2\pi i s}{\mu}} -  e^{\frac{2\pi i l}{\mu}}}{2}, ~ 
\frac{e^{\frac{2\pi i s}{\mu}} +  e^{\frac{2\pi i l}{\mu}}}{2}\right),
$$
where $s,l\in\{1,\ldots,\mu\}$ and $s\neq l$.
\end{itemize}
Formula (\ref{ParSumCoef}) with $d = 1$ gives
\begin{equation}
\label{ParSumCoef2} 
S_{[k]} = h_{2(k-1),k} = [t^{2(k-1)}]\frac{t^k(a+bt)^k}{(1-t)(1-t^2)},~~ k\geq 1.
\end{equation}
We have the following formula for such partial sums.
\begin{prop}
\label{prop9}
$$
S_{[k]} = \sum\limits_{n=0}^{k-2}{k\choose n} \left\lfloor \frac{k - n}{2} \right\rfloor b^na^{k-n}
$$
\end{prop}
\begin{proof}
Using the f.p.s.
$$
\frac{1}{(1 - t)(1 - t^2)} = \sum\limits_{i=1}^{\infty} i(t^{2i - 2} + t^{2i - 1}) = 
\sum\limits_{j=0}^{\infty} \left\lfloor 1 + \frac{j}{2}\right\rfloor t^j,
$$
and the main properties of the ``{\sl coefficient of }" operator 
(see \S 2.2 of \cite{Shapiro2}) we obtain 
$$
S_{[k]} = [t^{k - 2}]\frac{(a+bt)^k}{(1 - t)(1 - t^2)}= \sum\limits_{n=0}^{k-2}
[t^n] (a+bt)^k [t^{k-2-n}]\frac{1}{(1 - t)(1 - t^2)}  = 
$$
$$
\sum\limits_{n=0}^{k-2}{k\choose n}b^na^{k-n}\left\lfloor \frac{k - n}{2} \right\rfloor.
$$
\end{proof}
To see when these partial sums are eventually periodic, let us use Lemma \ref{lemma4} 
together with the rule of diagonalization (see \S 2.4 of \cite{Shapiro2}) saying 
\begin{equation}
\label{G6}
{\cal G} \bigl([t^n]F(t)\phi(t)^n\bigr) = \left[ \left. \frac{F(w)}{1 - t \phi'(t)} ~ 
\right| ~ w = t\phi(w) \right],
\end{equation}
where $\cal G$ denotes the {\sl generating function operator}.

\begin{theorem}
\label{theorem10}
Let $p(t) = a + bt \in\comp[t]$ be a polynomial, such that $a\cdot b\neq 0$. Then 
$V_{p(t)}$ has period $\mu$ and the sequence of column partial 
sums $\{S_{[k]}\}_{k\geq 1}$ is eventually periodic if and only if 
$$
-a = b = \frac{1}{2}e^{\frac{2\pi i s}{\mu}},~s\in\{1,\ldots,\mu\}.
$$
\end{theorem}
\begin{proof}
Suppose $V_{p(t)}$ has period $\mu$. To match formulas (\ref{ParSumCoef2}) and  (\ref{G6}), choose
$$
\phi(t) = a + bt ~~ \mbox{and} ~~ F(t) = \frac{(a + bt)^2}{(1 - t) (1-t^2)},
$$
then $w = ta/(1 - tb)$ and simple computations show that the generating 
function for the sequence $\{S_{[k]}\}_{k\geq 2}$ will be
\begin{equation}
\label{GFlinear}
\frac{F(w)}{1-tb} = \frac{a^2}{\bigl(1+ (a-b)t\bigr)\bigl(1 - (a+b)t\bigl)^2}.
\end{equation}
Now, since the numerator of ${\cal G}(S_{[k]})$ is a constant, 
Lemma \ref{lemma4} implies that the sequence $\{S_{[k]}\}_{k\geq1}$ can not be eventually periodic 
if the denominator of ${\cal G}(S_{[k]})$ has repeated roots, i.e. if $a + b \neq 0$. Thus $a = -b$, 
and we are in the Case (2). On the other hand, if
$$
- a = b = \frac{1}{2}e^{\frac{2\pi i s}{\mu}},~s\in\{1,\ldots,\mu\},
$$
we have
$$
V_{p(t)} = a\begin{pmatrix}
-1 & \phantom{-}1\\
\phantom{-}1 & -1\\
\end{pmatrix} ~~ \mbox{with} ~~ V^{n+1}_{p(t)} = (-2a)^nV_{p(t)},
$$
so $V_{p(t)}$ has period $\mu$, since $(-2a)^{\mu} = 1$. 
In this case (\ref{GFlinear}) simplifies to
$$
\frac{F(w)}{1-tb} = \frac{a^2}{1 + 2at} = a^2\bigl(1 + (-2a)t + (-2a)^2t^2 + 
(-2a)^3t^3 + \cdots \bigr),
$$
and we see that $\{S_{[k]}\}_{k\geq 1}$ is eventually periodic. Its general term  
is given by the explicit formula
\begin{equation}
\label{EQ-1}
S_{[k]} = a^2(-2a)^{k-2},~~\mbox{when} ~~k\geq 2.
\end{equation}
\end{proof}

As a simple application of this theorem and proposition above, we can obtain an elegant formula 
for the binomial transform of the sequence 
$$
\left\{\left\lfloor \frac{k - n}{2} \right\rfloor\right\}_{0\leq n \leq k}, ~~~ \forall k\in\natu.
$$

\begin{cor}
For all $k\geq 2$, we have
$$
\sum\limits_{n=0}^{k}(-1)^n{k\choose n} \left\lfloor \frac{k - n}{2} \right\rfloor  = (-2)^{k-2}.
$$
\end{cor}
\begin{proof}
The equality follows directly from Proposition \ref{prop9} and formula (\ref{EQ-1}), if 
we take $a = -b = 1$, since $\left\lfloor \frac{k - n}{2} \right\rfloor = 0$ when $n \in \{k-1,k\}$.
\end{proof}

%%%%%%%%%%%%%%%%%%%%%%%%%%%%%%%%%
%%%%%%%%%%%%%%%%%%%%%%%%%%%%%%%%%
%%%%%%%%                   			 %%%%%%%%%%%%
%%%%%%%%    Quadratic Polynomials   %%%%%%%%%%%%
%%%%%%%%                       			 %%%%%%%%%%%%
%%%%%%%%%%%%%%%%%%%%%%%%%%%%%%%%%
%%%%%%%%%%%%%%%%%%%%%%%%%%%%%%%%%

\subsection{Quadratic polynomials}

Take now a quadratic polynomial $p(t) = a + bt + ct^2$. Then 
$$
V_{p(t)} = \begin{pmatrix}
c & b & a\\
a & c & b\\
b & a & c\\
\end{pmatrix} ~~ \mbox{and} ~~ 
F^{-1} V_{p(t)} F = D = \begin{pmatrix}
\lambda_1 & 0 & 0\\
0 & \lambda_2 & 0\\
0 & 0 & \lambda_3\\
\end{pmatrix},
$$
where 
$$
F = \frac{1}{\sqrt{3}}\begin{pmatrix}
1 & 1 & 1\\
1 & \xi & \xi^2\\
1 & \xi^2 & \xi\\
\end{pmatrix},  ~~~~ F^{-1} = \frac{1}{\sqrt{3}}\begin{pmatrix}
1 & 1 & 1\\
1 & \xi^2 & \xi\\
1 & \xi & \xi^2\\
\end{pmatrix},
$$
and the eigenvalues $\lambda_i$ are 
$$
\lambda_1 = a + b + c, ~~~ \lambda_2 = a\xi^2 + b \xi + c, ~~~ \lambda_3 = a\xi + b\xi^2 + c.
$$
Therefore, if $V_{p(t)}$ has period $\mu$, then $V^{\mu+1}_{p(t)} = V_{p(t)}$, and
$$
FD^{\mu +1}F^{-1} = V_{p(t)}^{\mu + 1} = V_{p(t)} = FDF^{-1} \Longrightarrow D^{\mu+1} = D,
$$
which is equivalent to the system $\lambda_j^{\mu+1} = \lambda_j ,~ j\in\{1,2,3\}$, i.e. 
%(cf. \ref{PerCirMatLin})
\begin{equation}
\label{system3}
\left\{
\begin{array}{rcl}
(a+b+c)^{\mu+1} & = & a+b+c \\
(a\xi^2 + b \xi + c)^{\mu+1} & = & a\xi^2 + b \xi + c \\
(a\xi + b\xi^2 + c)^{\mu+1} & = & a\xi + b\xi^2 + c\\
\end{array}
\right.
\end{equation}
For a quadratic polynomial $p(t)$ with the matrix $V_{p(t)}$ of period $\mu$, we look at the sums 
$$
S_{3(k-1)} = \sum\limits_{i=0}^{3(k-1)}C_{i,k},
$$
which are the $3(k-1)$-st partial sums of the f.p.s.
$$
\sum\limits_{i=0}^{\infty}C_{i,k} t^i = \frac{\bigl(t(a+bt+ct^2)\bigr)^k}{1-t^3}.
$$
Let us follow the approach I used to prove Proposition 5 in \cite{KrylovRA1}. 
According to Theorems 1 and 2 of \cite{KrylovRA1}, the periodicity in $k$-th 
column of the $RAp$ begins with the $1 + 3(k-1)$-st term, and the repeating 
part equals the triple 
$$
\begin{pmatrix}
a_{k-1,0} \\ a_{k-1,1} \\ a_{k-1,2} \\
\end{pmatrix} = V^{k-1}_{p(t)} 
\begin{pmatrix}
a_0 \\ a_1 \\ a_2 \\
\end{pmatrix}, ~~ \mbox{where} ~~ a_0=a, a_1 = b, a_2 = c.
$$
Hence, if we use the polynomial $p_{k-1}(t) = a_{k-1,0} + a_{k-1,1}t + a_{k-1,2}t^2$, 
introduced in Definition \ref{definition3}, we can represent the infinite periodic 
tail of the $k$-th column by the generating function
$$ 
\left(t^{1 + 3(k-1)} \right)\frac{p_{k-1}(t)}{1-t^3}.
$$
Then we can write the partial sum $S_{3(k-1)}$ as
\begin{equation}
\label{MainPS}
S_{3(k-1)} = \left.\sum\limits_{i=0}^{3(k-1)}C_{i,k}t^i\right|_{t=1} = 
\lim\limits_{t\to 1}\left[\frac{\bigl(tp(t)\bigr)^k}{1-t^3} - 
\left(t^{1 + 3(k-1)} \right)\frac{p_{k-1}(t)}{1-t^3}\right].
\end{equation}
Assuming that $V_{p(t)}$ has period $\mu$ and the partial sums $\{S_{3(k-1)}\}_{k\geq 1}$ are 
eventually periodic with a period $h\in\natu$, we must have for each $k\in\{1,\ldots, h\}$ 
$$
S_{3(k-1)} = S_{3(k + nh - 1)}, ~\forall n\geq 0.
$$ 
In particular, for $n = \mu$ we have
$$
\lim\limits_{t\to 1} \left[\frac{\bigl(tp(t)\bigr)^k}{1-t^3} - 
\left(t^{1 + 3(k-1)} \right)\frac{p_{k-1}(t)}{1-t^3}\right]
$$
\begin{equation}
\label{PSums2}
= \lim\limits_{t\to 1}  \left[\frac{\bigl(tp(t)\bigr)^{k + h\mu}}{1-t^3} - 
\left(t^{1 + 3(k + h\mu -1)} \right)\frac{p_{k + h\mu -1}(t)}{1-t^3}\right].
\end{equation}
Since $1 - t^{3h\mu} = (1 - t^3)(1 + t^3 + t^6 + \cdots + t^{3(h\mu -1)})$, and 
$p_{k + h\mu -1}(t) = p_{k-1}(t)$ (because $\mu$ is the period of $V_{p(t)}$), 
we rewrite the difference of two rational functions
$$
\left(t^{1 + 3(k-1)} \right)\frac{p_{k-1}(t)}{1-t^3} - 
\left(t^{1 + 3(k + h\mu -1)} \right)\frac{p_{k + h\mu -1}(t)}{1-t^3} = 
t^{1 + 3(k-1)} \frac{p_{k-1}(t)}{1-t^3}\bigl( 1 - t^{3h\mu}\bigr) 
$$
as a polynomial 
\begin{equation}
\label{PSums3}
t^{1 + 3(k-1)}\cdot (1 + t^3 + t^6 + \cdots + t^{3(h\mu -1)})\cdot p_{k-1}(t).
\end{equation}
Moreover, since each bracket in (\ref{PSums2}), as well as (\ref{PSums3}), 
has a limit when $t \to 1$, the difference 
$$
\frac{\bigl(tp(t)\bigr)^k}{1-t^3} - \frac{\bigl(tp(t)\bigr)^{k + h\mu}}{1-t^3}
$$
also has a limit when $t \to 1$, and we can rewrite the equality in (\ref{PSums2}) as 
$$
\lim\limits_{t\to1} \left[\frac{(tp(t))^k\bigl(1 - (tp(t))^{h\mu}\bigr)}{1-t^3}\right]  -  
\lim\limits_{t\to1} \left[ t^{1 + 3(k-1)}
\left(1 - t^{3h\mu}\right)\frac{p_{k -1}(t)}{1-t^3}\right] = 0.
$$
Since $p_{k-1}(1) = (a+b+c)^k$ 
(recall Lemma \ref{lemma2}), the last identity is equivalent to
\begin{equation}
\label{EQ1}
(a+b+c)^k \lim\limits_{t\to1}\left[\frac{1 - (tp(t))^{h\mu}}{1-t^3}\right] - (a+b+c)^k(h\mu) = 0.
\end{equation}
We will consider the cases when $\lambda_1 = a + b + c\neq 0$ and $a + b + c = 0$ 
separately. Assume first that $a + b + c\neq 0$, then due to (\ref{system3}), 
$(a+b+c)^{\mu} = 1$. L'H\^opital's/Bernoulli's rule implies that
$$
\lim\limits_{t\to 1} \frac{1 - (tp(t))^{h\mu}}{1-t^3} = 
\lim\limits_{t\to 1} \frac{-h\mu (tp(t))^{h\mu - 1}\bigl(p(t) + tp'(t)\bigr)}{-3 t^2} = 
h\mu\frac{(a + 2b + 3c)}{3(a+b+c)},
$$
and we deduce from (\ref{EQ1}) that 
$$
a + 2b + 3c = 3a + 3b + 3c \Longleftrightarrow b = -2a
$$
(cf. proof of Proposition \ref{PropD}).
Hence we can eliminate $b$ from (\ref{system3}) to obtain
\begin{equation}
\label{system3b}
\left\{
\begin{array}{rcl}
(c - a)^{\mu} & = & 1\\
(a\xi^2 - 2a \xi + c)^{\mu+1} & = & a\xi^2 - 2a \xi + c \\
(a\xi - 2a \xi^2 + c)^{\mu+1} & = & a\xi - 2a\xi^2 + c\\
\end{array}
\right.
\end{equation}
If we assume for a moment that $\lambda_2 = a\xi^2+ b\xi + c = a\xi^2 - 2a\xi + c = 0$, then 
$$
1 = (c - a)^{\mu} = (2a\xi - a\xi^2 - a)^{\mu} = (3a\xi)^{\mu},
$$
and $\lambda_3 = a\xi - 2a \xi^2 + c = 3a\xi - 3a\xi^2 = 3a\xi( 1 - \xi)\neq 0$. Then 
the third equation in (\ref{system3b}) implies that $(3a\xi(1 - \xi))^{\mu} = 1 \Rightarrow 
(1 - \xi)^{\mu} = 1$, which is impossible since the norm of $1 - \xi$ equals $\sqrt{3}$. 
By a similar argument, $\lambda_3 = a\xi - 2a\xi^2 + c \neq 0$. Therefore, if 
$\lambda_1 \neq 0$, then none of the eigenvalues is zero, and (\ref{system3}) is equivalent to 
\begin{equation}
\label{system3c}
\left\{
\begin{array}{rcl}
(c - a)^{\mu} & = & 1\\
(a\xi^2 - 2a \xi + c)^{\mu} & = & 1 \\
(a\xi - 2a \xi^2 + c)^{\mu} & = & 1\\
\end{array}
\right. \Leftrightarrow 
\left\{
\begin{array}{rcl}
(c - a)^{\mu} & = & 1\\
(c - a - 3a\xi)^{\mu} & = & 1 \\
(c - a - 3a\xi^2)^{\mu} & = & 1\\
\end{array}
\right.
\end{equation}
Since $|\xi| = 1$, (\ref{system3c}) implies 
$$
\left\{
\begin{array}{rcl}
|(c - a)\xi | & = & 1\\
|(c - a)\xi^2 - 3a| & = & 1 \\
|(c - a)\xi - 3a| & = & 1. \\
\end{array}
\right.
$$
If we denote $(c-a)\xi$ and $3a$ by $w$ and $z$ correspondingly, we can 
apply Lemma \ref{lemma5} to see that $3a$ lies on the line through the origin, 
which bisects the angle $2\pi/3$ between the vectors $(c-a)\xi$ and $(c-a)\xi^2$. 
But $(c - a)\xi + (c - a)\xi^2 = a - c$ is on the same line, therefore there 
exists a nonzero scalar $r \in \real$ such that $c = ra$ and $c - a = a(r - 1)$. 
Furthermore, since $\lambda_1 = c - a \neq 0 \Rightarrow r -1 \neq 0$,
the equalities $(c - a)^{\mu} = 1 = (c - a - 3a\xi)^{\mu}$ imply  
$$
\left( 1 - \frac{3a\xi}{c - a}\right)^{\mu} = \left( 1 - \frac{3\xi}{r - 1} \right)^{\mu} = 1,
$$
so the norm of $1 - 3\xi/(r - 1)$ equals 1. Since $\bar{\xi} = \xi^2$ and $1 = -\xi- \xi^2$, 
$$
1 = \left( 1 - \frac{3}{r - 1} \xi\right)\left( 1 - \frac{3}{r - 1} \xi^2 \right) = 1 + 
\frac{3}{r - 1} + \frac{9}{(r - 1)^2}.
$$
The last equation has the only solution $r = -2$. Thus, we obtain $c = -2a$, and 
$p(t) = a(1 - 2t - 2t^2)$. If we rewrite system (\ref{system3}) in the form
\begin{equation}
\label{system3d}
\left\{
\begin{array}{rcl}
(- 3a)^{\mu} & = & 1\\
(- 3a(1 + \xi))^{\mu} & = & 1 \\
(- 3a(1+\xi^2))^{\mu} & = & 1,
\end{array}
\right.
\end{equation}
we deduce that $(-\xi^2)^{\mu} = (-\xi)^{\mu} = 1$, so $\mu$ must be 
divisible by 6. Let us summarize the case when $\lambda_1 \neq 0$ as follows.

{\sl 
If $V_{p(t)}$ has a finite period $\mu\in \natu$, such that the partial sums 
$\{S_{3(k-1)}\}_{k\geq 1}$ are eventually periodic and $p(1) = a + b + c \neq 0$, 
we must have 
$$
p(t) = a(1 - 2t - 2t^2), ~~~ \mbox{where} ~~~ (3a)^{\mu} = 1 ~\mbox{and} ~ 6~|~\mu. 
$$
}

On the other hand, if $p(t) = a(1 - 2t - 2t^2)$,  where $(3a)^{6m} = 1$, it is 
easy to see that $V^6_{p(t)}$ will have $(3a)^6$ on the leading diagonal, and 
zero elsewhere. Therefore $V_{p(t)} $ has the order (and the period) $6m$, and 
$p_{k - 1}(t) = p_{k + 6ms - 1}(t)$ for all $k\in\{1,\ldots, 6m\}$ and $s\in\natu$. 
Since for this polynomial $p(t)$ we have $a + 2b + 3c = 3(a + b + c) = -9a$, replacing 
$h\mu$ by $6ms$ in (\ref{PSums2}) - (\ref{EQ1}), one can use a similar  
argument as above, to show that $(\ref{PSums2})$ holds true, and hence the 
partial sums $\{S_{3(k-1)}\}_{k\geq 1}$ are eventually periodic with period $6m$.
Thus we proved

\begin{prop}
\label{propA}
Let $p(t) = a + bt + ct^2$ be a polynomial with $a\cdot c\cdot p(1) \neq 0$. Then 
$V_{p(t)}$ has a finite period $\mu\in \natu$ and the partial sums 
$\{S_{3(k-1)}\}_{k\geq 1}$ are eventually periodic if and only if 
$$
p(t) = a(1 - 2t - 2t^2), ~~~6~|~\mu, ~~ \mbox{and} ~~ (3a)^{\mu} = 1. 
$$
\end{prop}

Consider now the case when $\lambda_1 = a+b+c=0$, and 
(\ref{EQ1}) will be trivially true. If we also have $\lambda_2 = a\xi^2 + b\xi + c = 0$, 
then it follows from $1 + \xi + \xi^2 =0$ that 
$$
b = a\xi^2 ~~ \mbox{and} ~~ c = a\xi.
$$
Since $c\neq 0$, it is impossible for all three eigenvalues to be zero, and we 
deduce from (\ref{system3}) that when $\lambda_1 = \lambda_2 =0$, we have 
$$
(3a\xi)^{\mu} = 1,~~ \mbox{and} ~~ p(t) = a(1 + \xi^2 t + \xi t^2) = 
a(t - 1)(t - \xi^2)\xi = a(t - 1)(\xi t - 1).
$$
Using $1 - t^3 = (1-t)(1-\xi t)(1 - \xi^2 t)$ 
rewrite (\ref{ParSumCoef}) for $S_{[k]},\forall k\geq 2$, as
$$
S_{[k]} = [t^{2k-3}]\frac{\bigl(a(t - 1)(\xi t - 1)\bigr)^k}{(1-t)^2(1+t+t^2)} = 
a^{k} [t^{2k-3}]\frac{(1 - t)^{k-2}(1 - \xi t)^{k-1}}{1 - \xi^2 t},
$$
and apply formula (\ref{EQ3}) from Lemma \ref{lemma7} to obtain 
\begin{equation}
\label{Pr13A}
S_{[k]} = a^k3^{k-2}(\xi - 1)\xi^{k + 3 - 1} = (3a\xi)^k\frac{(\xi - 1)}{9\xi},~k\geq 2.
\end{equation}
Since $3a\xi$ is a root of unity, this sequence is clearly periodic (cf. (\ref{EQ-1})). 

Similarly, when $\lambda_1 = \lambda_3 =0$ we have $a = b\xi^2,~c = b\xi$, with 
$$
(3b\xi)^{\mu} = 1,~ p(t) = b(\xi^2  + t + \xi t^2) = 
b(t - 1)(t - \xi)\xi = b(t - 1)(\xi t - \xi^2).
$$
In this case the formulas (\ref{ParSumCoef}) and (\ref{EQ3B}) also 
give us a periodic sequence 
$$
S_{[k]} = [t^{2k-3}]\frac{\bigl(b(t - 1)(t - \xi)\xi)^k}{(1-t)^2(1+t+t^2)} = 
(\xi b)^k [t^{2k-3}]\frac{(1 - t)^{k-2}(\xi - \xi^3 t)^{k}}{(1 - \xi t)(1 - \xi^2 t)} = 
$$
\begin{equation}
\label{Pr13B}
(\xi^2 b)^k [t^{2k-3}]\frac{(1 - t)^{k-2}(1 - \xi^2 t)^{k-1}}{1 - \xi t} 
= (3b\xi)^k \frac{(1 - \xi)}{9}.
\end{equation}

On the other hand, if we take one of the polynomials $a(1 + \xi^2 t + \xi t^2)$ or 
$a(\xi^2  + t + \xi t^2)$, satisfying for some $\mu\in\natu$ respectively 
$(3a\xi)^{\mu} = 1$ or $(3a\xi^2)^{\mu} = 1$, 
one can see that the circulant matrix $V_{p(t)}$ will be singular, but with a 
finite period. For example for $p(t) = a(1 + \xi^2 t + \xi t^2)$, we have 
$$
\left(\frac{1}{3a\xi}V_{p(t)}\right)^2 = \frac{1}{3a\xi}V_{p(t)} = 
\begin{pmatrix}
1/3 & \xi/3 & \xi^2/3 \\
\xi^2/3 & 1/3 & \xi/3 \\
\xi/3 & \xi^2/3 & 1/3
\end{pmatrix},
$$
and hence $V^{\mu + 1}_{p(t)} = V_{p(t)}$. The explicit formulas we proved above 
for $S_{[k]}$ show that the partial sums $\{S_{3(k-1)}\}_{k\geq 1}$ for such 
polynomials will be eventually periodic. Thus we proved

\begin{prop}
\label{propB}
Let $p(t) = a + bt + ct^2$ be a polynomial such that $a\cdot c\neq 0$, $a+b+c=0$, and 
only one of the eigenvalues of $V_{p(t)}$ is different from zero. Then 
$V_{p(t)}$ has a finite period $\mu\in \natu$ and the partial sums 
$\{S_{3(k-1)}\}_{k\geq 1}$ are eventually periodic if and only if 
$$
p(t) = a(1 + \xi^2 t + \xi t^2), ~~ (3a\xi)^{\mu} = 1 ~~ \mbox{and} ~~ 
S_{[k]} = (3a\xi)^k\frac{(\xi - 1)}{9\xi},
$$
or 
$$
p(t) = a(\xi^2 + t + \xi t^2), ~~ (3a\xi^2)^{\mu} = 1 ~~ \mbox{and} ~~ 
S_{[k]} = (3a\xi^2)^k\frac{(1 - \xi)}{9},
$$
where $k\geq 2$.
\end{prop}

Let us now discuss the final case when $\lambda_1 = a+b+c=0$ and 
$\lambda_2\cdot \lambda_3\neq 0$. Then $c = - a - b$, and we can 
rewrite the last two equations of (\ref{system3}) as 
\begin{equation}
\label{EQ7}
\left\{
\begin{array}{lcl}
(a\xi^2 + b\xi - a - b)^{\mu} & = & 1\\ 
(a\xi + b\xi^2 - a - b)^{\mu} & = & 1\\
\end{array}
\right.
\Leftrightarrow
\left\{
\begin{array}{lcl}
(\xi - 1)^{\mu}(b - \xi^2 a)^{\mu} & = & 1\\ 
(\xi - 1)^{\mu}(a - \xi^2 b)^{\mu} & = & 1\\
\end{array}
\right.
\end{equation}
Dividing the first equation by the second one produces
\begin{equation}
\label{EQ8}
1 = \left(\frac{b - \xi^2 a}{a - \xi^2 b}\right)^{\mu} = 
\left(\frac{\frac{b}{a} - \xi^2}{1 - \xi^2 \frac{b}{a}}\right)^{\mu} = 
\left(\frac{-\xi\left(\frac{b}{a} - \xi^2\right)}{\left(\frac{b}{a} - \xi\right)}\right)^{\mu},
\end{equation}
which implies that the complex numbers 
\begin{equation}
\label{EQ8B}
\frac{b}{a} - \xi^2 ~~~ \mbox{and} ~~~ \frac{b}{a} - \xi
\end{equation}
have equal norms. Assuming for a moment that $b = 0$,  
deduce from (\ref{EQ7}) that
$$
\bigl(a(\xi^2 - 1)\bigr)^{\mu} =1 ~ \mbox{and} ~  \bigl(a(\xi - 1)\bigr)^{\mu} =1
\Rightarrow (\xi + 1)^{\mu} = (e^{2\pi i/6})^{\mu} = 1.
$$
Hence $6 ~ | ~ \mu$ and $p(t) = a(1 - t^2)$, where $\bigl(a(\xi - 1)\bigr)^{\mu} =1$. 
Formula (\ref{ParSumCoef}) for such $p(t)$ gives
$$
S_{[k]} = [t^{3k-3}]\frac{\bigl(ta(1-t^2)\bigr)^k}{(1-t)^2(1+t+t^2)} = 
a^k[t^{2k-3}]\frac{(1-t^2)^{k - 2}(1+t)^2}{(1+t+t^2)},
$$
and (\ref{EQ9}) from Lemma \ref{lemma6} implies that $\{S_{[k]}\}_{k\geq 1}$ 
is given by the formula
\begin{equation}
\label{Sk1}
\{S_{[k]}\}_{k\geq 1} = 
\left\{\frac{2}{3\sqrt{3}}(-a\sqrt{3})^k \sin\left(\frac{(k-4)\pi}{6}\right)\right\}_{k\geq 1}.
\end{equation}
Since $|a\sqrt{3}| = |a(\xi - 1)| = 1$, this sequence is clearly eventually periodic.
On the other hand, if we start with a polynomial $p(t) = a(1 - t^2)$, such 
that $6 ~ | ~ \mu$ and $\bigl(a(\xi -1)\bigr)^{\mu} =1$, 
then again, using (\ref{ParSumCoef}) with (\ref{EQ9}), 
we obtain the same formula (\ref{Sk1}), so 
the partial sums form a periodic sequence.

It is interesting to note that in this case, we could use (\ref{EQ10}), (\ref{G6}), 
and Lemma \ref{lemma4} to prove the periodicity of $\{S_{[k]}\}_{k\geq 1}$ 
without getting an explicit formula like (\ref{Sk1}). Indeed, we saw in the 
proof of Lemma \ref{lemma6} (recall the left-hand side of (\ref{EQ10})) that 
$\{S_{[k]}\}_{k\geq 1} =$ 
$$
a^k[x^{k-1}]\sum\limits_{i=0}^{k-2} (-1)^i{k-2\choose i}x^i(x - 
x^3+x^4 - x^6 + x^7 - x^9 + \cdots) 
$$
$$
= a^k[x^{k-1}]\sum\limits_{i=0}^{k-2} (-1)^i{k-2\choose i}x^i\frac{x(1 - x^2)}{1 - x^3} = 
a^k[x^{k-1}]\frac{x(1 - x^2)(1 - x)^{k-2}}{1 - x^3}
$$
\begin{equation}
\label{PSeq11}
= [x^{k-2}]\left(\frac{a^2(1 + x)(a - ax)^{k-2}}{1 + x + x^2}\right).
\end{equation}
To write the generating function 
\begin{equation}
\label{EQ12}
{\cal G} \left([x^{k-2}]\left(\frac{a^2(1 + x)(a - ax)^{k-2}}{1+x+x^2}\right)\right) = 
{\cal G} \left([x^{k-2}]\bigl(F(x)\phi^{k-2}(x)\bigr)\right)
\end{equation}
let
$$
\phi(x) = a - ax, ~ F(x) = \frac{a^2(1 + x)}{1 + x + x^2} ~~ \mbox{and} ~~ w = \frac{ax}{1 + ax}.
$$
Using the rule of diagonalization (\ref{G6}) we obtain the following formula for the 
G.F. ${\cal G}$ in (\ref{EQ12})
$$
{\cal G} = \frac{F(w(x))}{1-x\phi'(x)} = \frac{a^2(1 + 2ax)}{1 + 3ax + 3a^2x^2}.
$$
Furthermore, since 
$$
1 + 27a^6 x^6 = 1 - (ia\sqrt{3}x)^6 = (1-3ax+6a^2x^2-9a^3x^3+9a^4x^4)\cdot (1 + 3ax + 3a^2x^2)
$$
we can write ${\cal G} $ as the fraction
\begin{equation}
\label{Sk12}
\frac{a^2(1 + 2ax)(1-3ax+6a^2x^2-9a^3x^3+9a^4x^4)}{1 - (ia\sqrt{3}x)^6} = 
\frac{Q(x)}{1 - (ia\sqrt{3}x)^6},
\end{equation}
with some $Q(x)\in\comp[x]$. Since the norm $|ia\sqrt{3}| = 1$, Lemma 
\ref{lemma4} implies that ${\cal G}$ is an eventually periodic sequence.

~

Return now to (\ref{EQ8B}) and assume that $b\neq 0$. Thus, we have two 
nonzero complex numbers 
$$
\frac{b}{a} - \xi^2 ~~~ \mbox{and} ~~~ \frac{b}{a} - \xi
$$
of equal norms. Applying Lemma \ref{lemma5} to the numbers $z = b/a$ and 
$w = \xi$, we see that $z\in\real\setminus\{0\}$ and $b = za$. Denote this $z$ by $r$, 
then $c = -a(1 + r)$ and $p(t) = a\bigl(1 + rt - (1 + r)t^2\bigr) = 
a(1-t)(1 + (r+1)t)$, and (\ref{EQ7}) is equivalent to 
\begin{equation}
\label{EQ7b}
\left\{
\begin{array}{lcl}
\bigl((\xi - 1)a(r - \xi^2)\bigr)^{\mu} & = & 1\\ 
\bigl((\xi - 1)a(1 - \xi^2r)\bigr)^{\mu} & = & 1,\\ 
\end{array}
\right.
\end{equation}
which implies $(r - \xi^2)^{\mu} = (1 - \xi^2 r)^{\mu} \Longleftrightarrow 
(\xi r - 1)^{\mu} = (\xi  - r)^{\mu}$. According to Lemma \ref{lemma8}, if $r$ satisfies such 
equation, it will be a root of a degree $\mu$ polynomial with 
rational coefficients $f_{\mu}(x)\in\rati[x]$. Hence, we have at most $\mu$ such roots, and 
obviously, $r = 1$ will be a root for any $\mu\in\natu$.

To show that the sequence $\{S_{[k]}\}_{k\geq 1}$ is eventually periodic for such 
polynomials $p(t)= a(1-t)(1 + (r+1)t)$, let us use formula (\ref{MainPS})
$$
S_{3(k-1)} = \lim\limits_{t\to 1} \left[\frac{\bigl(tp(t)\bigr)^k}{1-t^3} - 
\left(t^{1 + 3(k-1)} \right)\frac{p_{k-1}(t)}{1-t^3}\right],
$$
which is greatly simplified in this case. Indeed, 
we proved in Lemma \ref{lemma2} that $p_m(1)=\bigl(a(1+r -(1+r))\bigr)^{m+1}= 0$,
therefore $(1 - t)~|~p_m(t)$ and we can write $p_{m}(t) = (1-t)\beta_{m}(t)$, 
for some $\beta_{m}(t)\in\comp[t]$. Then $S_{3(k-1)}$ 
$$
= \lim\limits_{t\to 1} \left[\frac{t(1-t)(1+(r+1)t)\bigl(tp(t)\bigr)^{k-1}}{(1-t)(1+t+t^2)} - 
\left(t^{1 + 3(k-1)} \right)\frac{(1-t)\beta_{k-1}(t)}{(1-t)(1+t+t^2)}\right]
$$
\begin{equation}
\label{MPS2}
= \frac{t^k}{1+t+t^2}\Biggl[\bigl(1 + (r+1)t\bigr)p^{k-1}(t) - t^{2(k-1)}\beta_{k-1}(t)\Biggr]_{t\to1} 
= \frac{-\beta_{k-1}(1)}{3},
\end{equation}
since $p^{k-1} (1) = 0$ when $k\geq 2$. Furthermore, since $V_{p(t)}$ is periodic 
with period $\mu$, we have $p_{k-1}(t) = p_{k-1 + n\mu}(t)$ for all $n\in\natu_0$, that is 
$$
(1 - t)\beta_{k-1}(t) = (1 - t)\beta_{k-1 + n\mu}(t).
$$
Hence $\beta_{k-1}(t) = \beta_{k-1 + n\mu}(t)$ for all $n\in\natu_0,$
and (\ref{MPS2}) implies that $S_{3(k-1)} = S_{3(k-1 + n\mu)}$, that is 
$\{S_{[k]}\}_{k\geq 1}$ is eventually periodic.

Notice that this argument gives a third proof of the periodicity of the partial 
sums $\{S_{[k]}\}_{k\geq 1}$ for the polynomials $p(t) = a (1 - t^2)$ 
(c.f. (\ref{Sk1}) and (\ref{Sk12})).

The converse statement is also clear, that is if 
\begin{equation}
\label{r_zero}
p(t) = a\bigl(1 + rt - (1 + r)t^2\bigr), 
\end{equation}
where $a$ and $r\in \real$ satisfy (\ref{EQ7b}), then $V_{p(t)}$ has 
period $\mu$ and the sequence $\{S_{[k]}\}_{k\geq 1}$ will be eventually periodic. 
If $r = 0$ in (\ref{r_zero}), $p(t) = a(1 - t^2)$ and (\ref{EQ7b}) implies that 
$\bigr((\xi - 1)a\bigl)^{\mu} = 1$ and $(-\xi^2)^{\mu} = 1$. Hence $6~|~\mu$, 
so we don't need to separate the case when $b = 0$ from the others. 
Thus we proved 

\begin{prop}
\label{propC}
Let $p(t) = a + bt + ct^2$ be a polynomial such that $a\cdot c\neq 0$, $a+b+c=0$, and 
the other two eigenvalues of $V_{p(t)}$ are different from zero. Then 
$V_{p(t)}$ has a finite period $\mu$ and the partial sums 
$\{S_{3(k-1)}\}_{k\geq 1}$ are eventually periodic if and only if 
$$
p(t) = a\bigl(1 + rt - (1 + r)t^2\bigr) = a(1 - t)\bigl(1 + (r+1)t\bigr),
$$
where $a$ and $r$ satisfy the system (\ref{EQ7b}), in which case $r \in \real $ 
is a root of a polynomial $f_{\mu}(t)$ of degree $\mu$ with rational coefficients.
\end{prop}

\noindent As for the polynomial $f_{\mu}(t)$, we proved in Lemma \ref{lemma8} 
that it is given by the formula
$$
f_{\mu}(t) = \sum\limits_{i=0}^{\mu} (-t)^i{\mu\choose i} \cdot 
\left( \frac{\left(\xi^i + \bar{\xi}^i\right) - (-1)^{\mu}\left(\xi^{\mu - i} +\bar{\xi}^{\mu - i}\right)}
{2 - (-1)^{\mu}\bigl(\xi^{\mu} + \bar{\xi}^{\mu}\bigr)}   \right),~~~ \mbox{if} ~~~ 6 \nmid \mu,
$$
and 
$$
f_{\mu}(t) = \sum\limits_{i=0}^{\mu} (-t)^i{\mu\choose i} \cdot 
\frac{\xi^i - \bar{\xi}^i}{\sqrt{-3}},~~~ \mbox{if} ~~~ 6 ~|~ \mu.
$$
Here is the table of the first few such polynomials including their roots.
\begin{center}
\begin{tabular}{|c|c|c|}
\hline
${\mu}$ & $f_{\mu}(t)$ & Roots $r$\\
\hline
\hline
1 & $1 - t$ & $\{1\}$\\
\hline
2 & $1 - t^2$ & $ \{-1, ~1\}$\\
\hline
3 & $1 + (3/2)t - (3/2)t^2 - t^3$ & $\{-2,  -1/2,  1\}$\\
\hline
4 & $1 + 4t - 4t^3 - t^4$ & $ \{-2 - \sqrt{3},  -1,  -2 + \sqrt{3},  1\}$\\
\hline
5 & $1 + 10t + 10t^2 - 10t^3 - 10 t^4 - t^5$ & $\{-8.74, -1.46,-0.68, -0.11,  1\}$\\
\hline
6 & $-6t - 15t^2 + 15t^4 + 6t^5$ & $\{-2, -1,  -1/2,  0,  1\}$\\
\hline
\end{tabular}
\end{center}
It is also interesting to look at the corresponding integer-coefficient array 
of these polynomials. The integer coefficients are obtained by multiplying 
through by the l.c.m. of all denominators, when needed. This array is not a 
Riordan array, since $\xi^{6k} - \bar{\xi}^{6k} = 0$ for all $k\in\natu_0$, 
so each $6k$-th row has zero on the main diagonal.
$$
\setcounter{MaxMatrixCols}{25}
\begin{matrix}
\phantom{-}1 & -1 \\
\phantom{-}1 & \phantom{-}0 & -1\\
\phantom{-}2 &  \phantom{-}3 & -3 & -2\\
\phantom{-}1 & \phantom{-1}4 & \phantom{-}0 & -4 & -1\\
\phantom{-}1 & \phantom{-1}10 & \phantom{-} 10 & -10 & -10 & -1\\
\phantom{-}0 & \phantom{-1}6 & \phantom{-}15 & \phantom{-} 0 & -15 & -6 & \phantom{-}0\\
\phantom{-}1 & -7 & -42 & -35 & \phantom{-}35 & \phantom{-}42 &  \phantom{-}7 & -1\\
\phantom{-}1 & \phantom{-}0 & -28 & -56 & \phantom{-}0 & \phantom{-}56 &  \phantom{-}28 & 
\phantom{-}0 & -1\\
\phantom{-}2 & \phantom{-}9 & -36 & -168 & -126 & \phantom{-}126 & \phantom{-}168 &  
\phantom{-}36 & -9 & -2\\
\end{matrix}
$$

%%%%%%%%%%%%%%%%%%%%%%%%%%%
%%%%%%%%%%%%%%%%%%%%%%%%%%%
%%%%%%%%                       %%%%%%%%%%%%
%%%%%%%%     Section 4   %%%%%%%%%%%%
%%%%%%%%                       %%%%%%%%%%%%
%%%%%%%%%%%%%%%%%%%%%%%%%%%
%%%%%%%%%%%%%%%%%%%%%%%%%%%

\section{Polynomials with degree $d \geq 2$}

In this section we generalize the case of $p(t) = a(1 - 2t - 2t^2)$ to polynomials of 
an arbitrary degree, and show that the corresponding column partial sums are 
eventually periodic. Let us fix an integer $d\geq 2$ and consider the polynomial
\begin{equation}
\label{Polynom_d}
p(t) = a\bigl((d - 1) - \sum\limits_{i = 1}^d{2t^i}\bigr), ~~~ 
\mbox{where} ~~~ \bigl((d+1)a\bigr)^{k} = 1 ~\mbox{for some} ~ k\in\natu.
\end{equation}

\begin{theorem}
If $a = 1/(d+1)$, the circulant matrix $V_{p(t)}$ is nonsingular and has order $d+1$ if $d$ 
is odd, or $2(d+1)$ if $d$ is even respectively.
\end{theorem}
\begin{proof}
It is well known (see \cite{Kra}) that the eigenvalues of the circulant matrix generated by a vector 
$\nu = \{\nu_0\ldots, \nu_d\}$ are
$$
\lambda_l = \nu_0 + \xi^l\nu_1 + \cdots + \xi^{ld}\nu_d, ~~ \mbox{for every} ~~ l\in\{1,\ldots, d + 1\}
$$
(recall that $\xi = e^{2\pi i/(d+1)}$ here). For our polynomial $p(t)$ we have 
$$
\nu = \frac{1}{d+1}\{-2,-2,\ldots, -2, d-1\},
$$
and it follows from the identity $0 = \xi^{l(d+1)} - 1 = (\xi^l - 1)(1 + \xi^l + \cdots + \xi^{ld})$ 
that for all $l\in\{1,\ldots,d\}$, 
$$
\xi^l\lambda_l = \frac{-2}{d + 1}\left(\xi^l + \xi^{2l} + \cdots + \xi^{dl}\right) + 
\frac{d - 1}{d + 1} = \frac{2}{d + 1} +  \frac{d - 1}{d + 1} = 1.
$$ 
As for $l = d+1$, it is easy to see that $\lambda_{d+1} = -1$. Thus, for 
our polynomial $p(t)$, all the eigenvalues are roots of unity
$$
\lambda_{d+1} = -1~~ \mbox{and} ~~ \lambda_l = e^{\frac{2\pi i(d+1 - l)}{d+1}}, ~ 
\forall l\in\{1,\ldots, d\}.
$$
Hence, if $d+1$ is even, $\lambda_l^{d+1} = 1$ for every index $l$, and the matrix $V_{p(t)}$ 
has order $d+1$. When $d+1$ is odd, we have $\lambda_{d+1}^{d+1} = -1$, and we have 
to take $V_{p(t)}$ to the power of $2(d+1)$ to get the identity matrix.
\end{proof}
\noindent Immediately from this theorem we obtain the following 
\begin{cor}
\label{cor17}
If $a(d + 1)$ is a $k$-th root of unity, the matrix $V_{p(t)}$ has order $\lcm{(k, d+1)}$ or $\lcm{(k,2(d+1))}$.
\end{cor}

Now we generalize the argument used in \S 3.2 to prove (\ref{EQ1}), to show that 
for the polynomial (\ref{Polynom_d}) the corresponding sequence of column partial 
sums will be eventually periodic.

\begin{prop}
\label{PropD}
Take any integer $d\geq 2$, and $a\in\comp$ s.t.
$\bigl((d+1)a\bigr)^{k} = 1$, where $k$ is the smallest such positive integer.
Then the sequence of column partial sums $\{S_{[k]}\}_{k\geq 1}$ of the 
polynomial $p(t) = a\bigl((d - 1) - \sum\limits_{i = 1}^d{2t^i}\bigr)$ is 
eventually periodic with the period 
$\lcm{(k, d + 1)}$ when $d$ is odd, or $\lcm{(k, 2(d + 1))}$ when $d$ is even.
\end{prop}
\begin{proof}
Since the partial sum of the $k$-th column equals 
$$
S_{(k-1)(d+1)} = \sum\limits_{i=0}^{(k-1)(d+1)} C_{i,k},
$$
(recall (\ref{PSum})) the analog of formula (\ref{MainPS}) for  
$p(t) = a_0 + a_1 t + \cdots + a_d t^d$ is 
\begin{equation}
\label{MainPS2}
\left.\sum\limits_{i=0}^{(k -1)(d + 1)}C_{i,k}t^i\right|_{t=1} = 
\lim\limits_{t\to 1}\left[\frac{\bigl(tp(t)\bigr)^k}{1-t^{d+1}} - 
\left(t^{1 + (k-1)(d+1)} \right)\frac{p_{k-1}(t)}{1 - t^{d+1}}\right],
\end{equation}
where $p_{k-1}(t) = \sum\limits_{i = 0}^d a_{k - 1,i}t^i$ is the polynomial  
introduced in Definition \ref{definition3}. It satisfies the periodicity 
$p_{k+h\mu -1}(t) = p_{k-1}(t)$, with $\mu$ being the period of $V_{p(t)}$ 
and $h$ an arbitrary nonnegative integer. To make the formulas look less 
cumbersome, we will use the following notations:
$$
A_{k}:= \frac{\bigl(tp(t)\bigr)^k}{1-t^{d+1}}, ~~~ 
B_{k}:= \left(t^{1 + (k-1)(d+1)} \right)\frac{p_{k-1}(t)}{1 - t^{d+1}},
$$
and
$$
\mu : = \left\{
\begin{array}{lcl}
\lcm{(k,d+1)} & \mbox{if} & 2\nmid d\\
\lcm{(k,2(d+1))} & \mbox{if} & 2 ~ | ~ d\\
\end{array}\right..
$$
Then, as follows from Corollary \ref{cor17} and (\ref{MainPS2}), it is enough to show that
\begin{equation}
\label{PSums2B}
\lim\limits_{t\to 1}(A_{k} - B_{k})  -  
\lim\limits_{t\to 1}(A_{k+h\mu} - B_{k+h\mu}) = 0,~\forall h\in\natu_0.
\end{equation}
Take any $h\in\natu_0$. Since 
$1 - t^{h\mu (d+1)} = (1 - t^{d+1})(1 + t^{d+1} + \cdots + t^{(d+1)(h\mu -1)})$ and 
$p_{k + h\mu -1}(t) = p_{k-1}(t)$, we have 
$$
B_{k+h\mu} - B_k = B_k(t^{h\mu (d+1)} - 1)
$$
$$
= -\left(t^{1 + (k-1)(d+1)} \right)p_{k-1}(t)(1 + t^{d+1} + \cdots + t^{(d+1)(h\mu -1)}).
$$
Applying $\bigl((d+1)a\bigr)^{k} = 1$ and Lemma \ref{lemma2} to our polynomial $p(t)$ with 
$$
p_{k-1}(1) = \left(\sum\limits_{i=0}^d a_i \right)^{k} = \bigl(a(d - 1 - 2d)\bigr)^{k} = 
\bigl(-a(d + 1)\bigr)^{k} = (-1)^{k},
$$
we obtain 
\begin{equation}
\label{MainLimit}
\lim\limits_{t\to 1} (B_{k+h\mu} - B_k) = (-1)^{k+1} h\mu.
\end{equation}
Since all limits in (\ref{PSums2B}) and (\ref{MainLimit}) exist, 
$\lim\limits_{t\to1} (A_k - A_{k+h\mu})$ also exists. Hence, using 
L'H\^opital's/Bernoulli's rule we have 
$$
\lim\limits_{t\to1} (A_k - A_{k+h\mu}) = \lim\limits_{t\to 1} 
\left[\frac{(tp(t))^k\bigl(1 - (tp(t))^{h\mu}\bigr)}{1 - t^{d+1}}\right] = 
(-1)^k\lim\limits_{t\to1}\left[\frac{1 - (tp(t))^{h\mu}}{1 - t^{d+1}}\right]
$$
$$
= (-1)^k \lim\limits_{t\to 1}\frac{-h\mu (tp(t))^{h\mu - 1}\bigl(p(t) + tp'(t)\bigr)}{-(d+1)t^d} 
= (-1)^k \frac{h\mu p(1)^{h\mu}}{(d+1)p(1)}\bigl(p(1) + p'(1)\bigr)
$$
$$
=\frac{(-1)^k h \mu}{-a(d + 1)^2}\bigl( -a(d+1) - ad(d+1)\bigr) = (-1)^kh\mu,
$$
since $p(1)^{\mu} = 1$ and $p'(1) = -2a(1 + 2 + 3 + \cdots + d) = - ad(d + 1)$. 
Thus, we proved that for an arbitrary $h\in\natu_0$,
$$
S_{(k-1)(d+1)} - S_{(k-1+h\mu)(d+1)} = \lim\limits_{t\to 1}(A_{k} - B_{k}) -  
\lim\limits_{t\to 1}(A_{k+h\mu} - B_{k+h\mu})
$$
$$
= \lim\limits_{t\to1} (A_k - A_{k+h\mu}) + \lim\limits_{t\to 1} (B_{k+h\mu} - B_k) = 
(-1)^kh\mu + (-1)^{k+1} h\mu = 0,
$$
which confirms the periodicity and finishes the proof.
\end{proof}

\begin{example}
\label{example1}
Take $d = 3$ and $k=5$. Then 
$$
a = 1/4 ~~~ \mbox{and} ~~~ a = e^{\frac{6\pi i}{5}}/4
$$
are particular solutions of $\bigl(a(d+1)\bigr)^k = 1$, and the corresponding 
polynomials are $p(t) = a(2 - 2t - 2t^2 - 2t^3)$. For $a = 1/4$ the sequence 
of column partial sums is periodic with period $4$ and the repeating part
$$
\{0,\frac{-1}{2},\frac{1}{2},0\}.
$$
It has period 4 because for $a = 1/4$, $k=1$ is the smallest positive integer satisfying 
$\bigl(a(d+1)\bigr)^k = 1$, and $\lcm{(1,4)} = 4$. When $a = e^{\frac{6\pi i}{5}}/4$, 
the period will be $20 = \lcm{(5,4)}$, and the repeating part equals 
$$
\{0, \frac{-\xi^2}{2}, \frac{-\xi^3}{2}, 0, 0, \frac{\xi}{2}, \frac{\xi^2}{2},  0, 0,  \frac{-1}{2}, 
\frac{-\xi}{2} , 0, 0,  \frac{-\xi^4}{2}, \frac{1}{2}, 0, 0, \frac{\xi^6}{2}, \frac{\xi^8}{2}, 0\},
$$
where $\xi = e^{\frac{2\pi i}{10}}$. For graphical representation of these two sequences, 
see the first two graphs in Figure \ref{fig6} in the next section.
\end{example}

%%%%%%%%%%%%%%%%%%%%%%%%%%%
%%%%%%%%%%%%%%%%%%%%%%%%%%%
%%%%%%%%                       %%%%%%%%%%%%
%%%%%%%%     Section 5   %%%%%%%%%%%%
%%%%%%%%                       %%%%%%%%%%%%
%%%%%%%%%%%%%%%%%%%%%%%%%%%
%%%%%%%%%%%%%%%%%%%%%%%%%%%

\section{Graphs of column partial sums}

Every eventually periodic sequence consists of finitely many values. When 
these values are complex numbers, it is natural to connect two points in $\comp$ 
by an edge, when the numbers correspond to the consecutive elements of 
the sequence. This way, each eventually periodic sequence $\{S_{[k]}\}_{k \geq  1}$ 
generates a finite graph in the plane. In this section I will show several examples 
of such graphs representing the column partial sums we discussed above. 
All graphs are produced using Wolfram Mathematica. 

\subsection{Polynomials $p(t) = a(1 - t)$}

\begin{figure}[hbt!] %  figure placement: here, top, bottom, or page
\centering
\includegraphics[width=134mm]{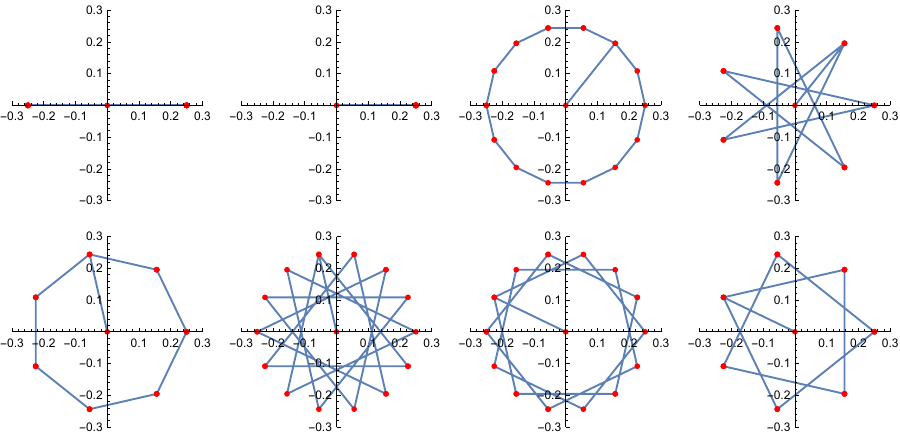} 
\caption{$p(t) = a(1 - t)$, where $-2a = e^{\frac{2\pi i s}{14}}$} 
\label{fig1}
\end{figure}

According to Theorem \ref{theorem10}, such a polynomial generates an 
eventually periodic sequence when
$$
a = -e^{\frac{2\pi i s}{\mu}}/2,~~s\in\{1,2,\ldots,\mu\},
$$
in which case the terms are given by
$$
S_{[k]} = a^2(-2a)^{k-2},~~~\mbox{where} ~~~ k\geq 2 ~\mbox{and}~ (-2a)^{\mu} = 1.
$$
Since $(-2a)$ is a root of unity, all vertices of the corresponding 
graph (except the first one, since $S_{[1]} = 0$) will be on a circle, subdividing it 
into congruent arcs. The vertices will be connected by an edge according 
to a choice of $s\in \{1,2,\ldots,\mu\}$. If $a$ solves the equation 
$(-2a)^{\mu} = 1$, its conjugate $\bar{a}$ is a solution too. If we change $a$ for its 
conjugate in $p(t)= a(1 - t)$, formula (\ref{EQ-1}) implies that 
each element of the sequence $\{S_{[k]}\}_{k\geq 1}$ also gets 
changed to its conjugate, and therefore the corresponding graphs will be 
symmetric about the $x$-axis. Hence, for our graphical presentation we can take 
only one element from each pair of the conjugate solutions of $(-2a)^{\mu} = 1$. 
Figure \ref{fig1} shows several examples of such graphs when $\mu = 14$. 
The first two graphs with three and two vertices correspond to $a = 1/2$ and 
$a = -1/2$, and lie entirely on the $x$-axis. 

\subsection{Polynomials $a(1 - 2t - 2t^2)$}

This is a case where we don't have an explicit formula for the sequence of 
column partial sums, but we can easily obtain several first terms using modern 
technology. Recall that here $6~|~\mu, ~~~ \mbox{and} ~~~ (3a)^{\mu} = 1$. 
If we take $\mu = 12$, we will have two real solutions and five pairs of complex 
conjugate numbers. Here is a list of seven different solutions of $(3a)^{12} = 1$ 
(excluding the conjugates)
$$
a\in\{\frac{-1}{3}, \frac{1}{3}, \frac{i}{3}, \frac{e^{\pi i / 6}}{3}, \frac{e^{\pi i / 3}}{3}, 
\frac{e^{2\pi i / 3}}{3},  \frac{e^{5\pi i /6}}{3}\},
$$
with the list of the corresponding graphs in Figure \ref{fig2} below.
\begin{figure}[hbt!]  %  figure placement: here, top, bottom, or page
\centering
\includegraphics[width=134mm]{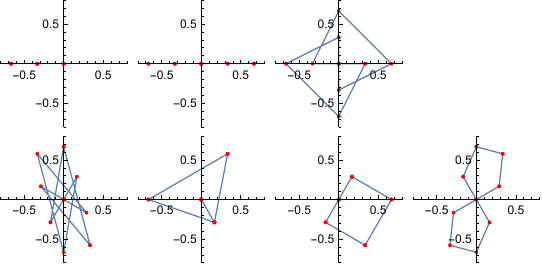} 
\caption{$p(t) = a(1 - 2t - 2t^2)$, where $(3a)^{12}=1$} 
\label{fig2}
\end{figure}

\noindent For example, the third graph corresponds to $a = i/3$, and represents the 
sequence $\{S_{[k]}\}_{k\geq 1}$ of period 12 with the repeating part
$$
\{0, \frac{1}{3}, \frac{-2i}{3}, \frac{-2}{3}, \frac{i}{3}, 0, 0,  \frac{-1}{3}, 
\frac{2i}{3},  \frac{2}{3}, \frac{-i}{3}, 0\}.
$$

\subsection{Polynomials $a(1 + \xi^2 t + \xi t^2)$ and $a(\xi^2 + t + \xi t^2)$}

\begin{figure}[hbt!] %  figure placement: here, top, bottom, or page
\centering
\includegraphics[width=134mm]{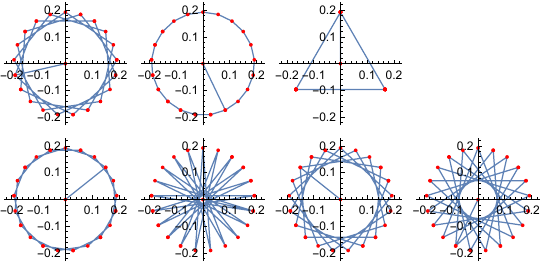} 
\caption{$p(t) = a(\xi^2 + t + \xi t^2)$, where $(3a\xi^2)^{7}=1$} 
\label{fig3}
\end{figure}

Here the terms of $\{S_{[k]}\}_{k\geq 1}$ are given explicitly by 
$$ 
S_{[k]} = (3a\xi)^k\frac{(\xi - 1)}{9\xi} ~~~ \mbox{with} ~~~ (3a\xi)^{\mu} = 1,
$$
and
$$
S_{[k]} = (3a\xi^2)^k\frac{(1 - \xi)}{9} ~~~ \mbox{with} ~~~ (3a\xi^2)^{\mu} = 1
$$
respectively (recall (\ref{Pr13A}) and (\ref{Pr13B})). Similar to the linear case, all 
vertices of each graph will be on a circle subdividing it into congruent arcs. 
For example, taking $\mu = 7$ we obtain seven complex solutions of 
$(3a\xi^2)^{\mu} = 1$ with the list of the 
corresponding graphs in Figure \ref{fig3}.

\noindent The third graph corresponds to $a = \xi/3$, and represents the 
sequence $\{S_{[k]}\}_{k\geq 1}$ of period 3 with the first four terms
$$
\{0, \frac{i}{3\sqrt{3}}, \frac{3 - i\sqrt{3}}{18}, \frac{-3 - i\sqrt{3}}{18}\}.
$$

\subsection{Polynomials $a(1 + rt - (1 + r)t^2)$, with $r\in\real$}

When $r = 0$, formula (\ref{Sk1}) implies that the vertices corresponding to the 
terms of $\{S_{[k]}\}_{k\geq 1}$ will be again on the rays starting at the origin 
and dividing the rotation angle ($2\pi$) into equal sub-angles, so this case 
is more or less similar to the ones we saw above. The difference will be in the 
distance from the origin, which will vary according to the function 
$\sin\left(\frac{(k-4)\pi}{6}\right)$, and so we focus on $r\neq 0$. Furthermore, 
it is straightforward to check that $a$ is a solution of the system 
(\ref{EQ7b}) if and only if $\bar{a}$ is a solution as well,
so here we also can take only one representative from each pair of 
the conjugate solutions. As we saw in Lemma \ref{lemma8}, when $6~|~\mu$ 
the formula for the polynomial $f_{\mu}(t)$ is different, so let us choose two values,
one of which is divisible by 6, e.g. $\mu = 12$ and $\mu = 14$. 

\begin{figure}[hbt!] %  figure placement: here, top, bottom, or page
\centering
\includegraphics[width=134mm]{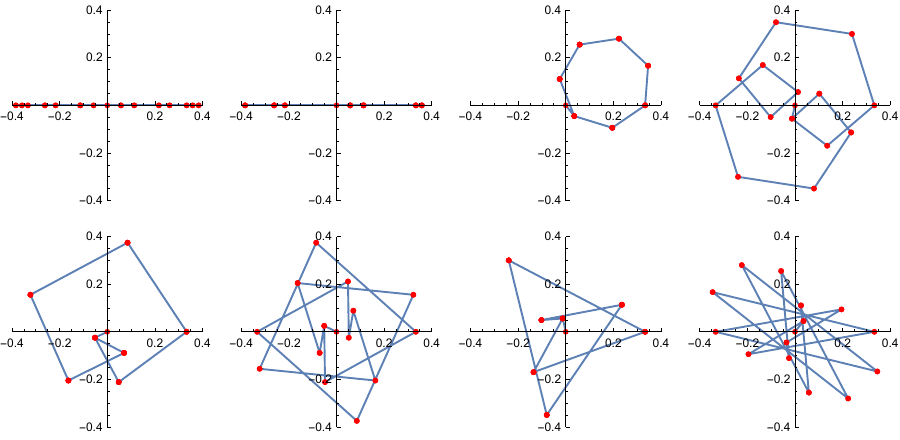} 
\caption{$p(t) = a(1 - t)(1 + 12.0563 t)$, with $\mu = 14$} 
\label{fig4} ~ \\ ~ \\ ~ 
\includegraphics[width=124mm]{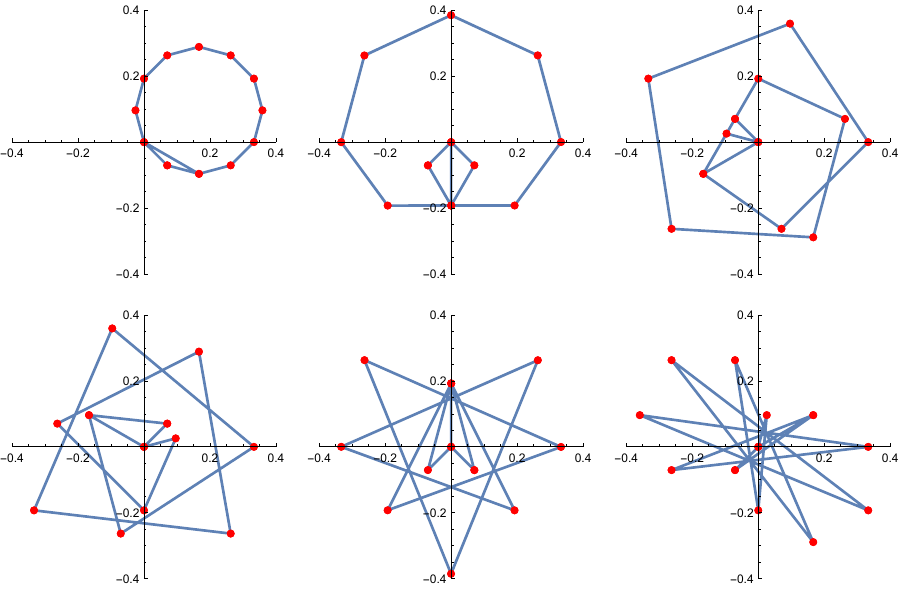} 
\caption{$p(t) = a(1 - t)(1 + 3.732 t)$, with $\mu = 12$} 
\label{fig5}
\end{figure}

If $\mu = 14$, $f_{14}(t)$ will have 14 real roots, with the 
largest one $r\approx 11.0563$. The list of the corresponding graphs is 
given in Figure \ref{fig4}.

If $\mu = 12$, $f_{12}(t)$ will have 11 roots
$$
\{-2, -1, \frac{-1}{2}, 0, 1, -2 - \sqrt{3}, \frac{-1 - \sqrt{3}}{2}, 
 1 - \sqrt{3}, -2 + \sqrt{3}, \frac{-1 + \sqrt{3}}{2}, 1 + \sqrt{3}\},
$$
with the 
largest one $r = 1 + \sqrt{3} \approx 2.732$. The list of the corresponding graphs is 
given in Figure \ref{fig5} below.

\subsection{Polynomials $a(2 - 2t - 2t^2 - 2t^3)$}

Here we also take only one element from each pair of the 
conjugate solutions of $(4a)^k = 1$. Recall that in Example \ref{example1}, we took 
$k = 5$ so for this case we will have three (essentially) different graphs. They are 
shown in Figure \ref{fig6} below. The first two graphs represent sequences 
with periods 4 and 20, which are generated by $a = 1/4$ and $a = e^{\frac{6\pi i}{5}}/4$, 
respectively (Example \ref{example1}).

\begin{figure}[hbt!] %  figure placement: here, top, bottom, or page
\centering
\includegraphics[width=133mm]{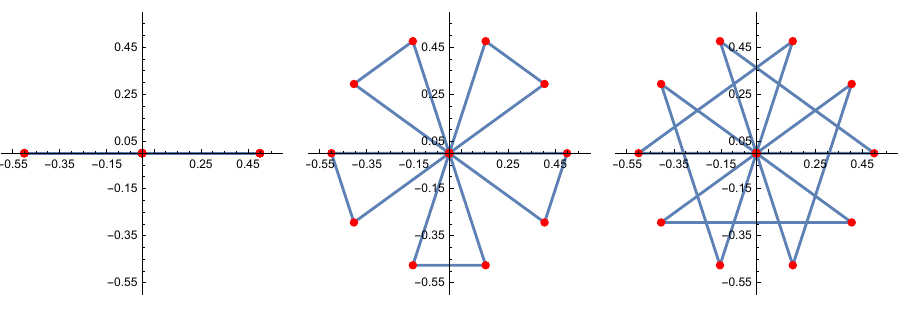} 
\caption{$p(t) = a(2 - 2t - 2t^2 - 2t^3)$, where $(4a)^5=1$} 
\label{fig6}
\end{figure}

\end{document}